\documentclass{article}

\usepackage[preprint]{neurips_2021}

\usepackage[utf8]{inputenc} 
\usepackage[T1]{fontenc}    
\usepackage{hyperref}       
\usepackage{url}            
\usepackage{booktabs}       
\usepackage{amsfonts}       
\usepackage{nicefrac}       
\usepackage{microtype}      
\usepackage{xcolor}         
\usepackage{bm, bbm}
\newtheorem{definition}{Definition}
\newtheorem{lemma}{Lemma}
\newtheorem{assumption}{Assumption}
\newtheorem{proposition}{Proposition}
\newtheorem{theorem}{Theorem}
\newtheorem{remark}{Remark}
\usepackage{amsmath,amssymb,mathtools}
\usepackage{multirow,rotating,makecell,longtable,threeparttable,booktabs,caption,lscape, graphicx, geometry, subfig}

\newcommand{\Fscr}{\ensuremath{\mathcal F}}

\newcommand{\Wscr}{\ensuremath{\mathcal W}}
\newcommand{\Xscr}{\ensuremath{\mathcal X}}
\newcommand{\Yscr}{\ensuremath{\mathcal Y}}
\newcommand{\Zscr}{\ensuremath{\mathcal Z}}
\def\E{\mathbb{E}}

\newcommand{\field}[1]{\ensuremath{\mathbb{#1}}}
\newcommand{\R}{\ensuremath{\field{R}}} 
\newcommand{\I}[1]{\ensuremath{\mathbb{I}_{\left\{#1\right\}}}} 

\DeclareMathOperator*{\argmin}{arg\,min}
\usepackage{float}
\newcommand{\cny}[1]{\textcolor{red}{NY: #1}}
\newcommand{\wty}[1]{\textcolor{orange}{TY: #1}}
\title{Distributionally Robust Prescriptive Analytics with Wasserstein Distance}

\author{%
  Tianyu Wang 
   \\
  School of Economics and Management\\
  Tsinghua University\\
  \texttt{wangtianyu6162@gmail.com} \\
  \And
  Ningyuan Chen\\
  Rotman School of Management\\
  University of Toronto\\
  \texttt{ningyuan.chen@utoronto.ca} \\
  \And
  Chun Wang\\
  School of Economics and Management\\
  Tsinghua University\\
  \texttt{wangch5@sem.tsinghua.edu.cn} \\
}

\begin{document}

\maketitle

\begin{abstract}
In prescriptive analytics, the decision-maker observes historical samples of $(X, Y)$, where $Y$ is the uncertain problem parameter and $X$ is the concurrent covariate, without knowing the joint distribution. Given an additional covariate observation $x$, the goal is to choose a decision $z$ conditional on this observation to minimize the cost $\mathbb{E}[c(z,Y)|X=x]$. This paper proposes a new distributionally robust approach under Wasserstein ambiguity sets, in which the nominal distribution of $Y|X=x$ is constructed based on the Nadaraya-Watson kernel estimator concerning the historical data.  We show that the nominal distribution converges to the actual conditional distribution under the Wasserstein distance.  We establish the out-of-sample guarantees and the computational tractability of the framework.  Through synthetic and empirical experiments about the newsvendor problem and portfolio optimization, we demonstrate the strong performance and practical value of the proposed framework.
\end{abstract}

\section{Introduction}
Most decision-making problems involve uncertainty. 
In an increasingly data-rich world, the decision-maker may rely on covariates or contextual information to predict the uncertain information (though imperfectly) when making such a decision.
This is the central question in predictive analytics, which has attracted the attention of many scholars.


\paragraph{Prescriptive Analytics}
Consider the following optimization problem:
\begin{equation}\label{eq:formulation}
    \min_{z \in \Zscr} \E[c(z, Y) | X = x].
\end{equation}
Here $z$ denotes the decision vector taken from the feasible region $\Zscr\subset \R^{d_z}$.
The objective function $c(z,Y)$ is known and depends on the random vector $Y\in \R^{d_y}$, which encodes the information of the environment when the decision has to be made.
However, $Y$ is typically unobserved and the decision-maker has access to a covariate $X\in \R^{d_x}$.
The joint distribution $(X,Y)$ is usually not known. 
Instead, in most applications, the decision maker can observe $n$ samples $\{(x_i, y_i)\}_{i = 1}^{n}$, presumably from the decisions made in the past.
The goal of prescriptive analytics is to choose a reasonable decision $z$ based on the current covariate $X=x$ and historical data $\{(x_i, y_i)\}_{i = 1}^{n}$.

We give some applications of the formulation in the table below.
Because of the practical relevance, the problem has been investigated in many recent papers \cite{elmachtoub2017smart,tulabandhula2013machine,bertsimas2020predictive,bertsimas2021data,notz2021prescriptive}. 
The novelty of this paper is the methodology we use to handle the problem, which is described next.
\begin{center}
   \begin{tabular}{c|ccc}
   \toprule
     Problem &$X$&$Y$&$  z$\\
     \midrule
     Newsvendor&Weather&Demand&Order quantity\\
     Personalized pricing&Consumer feature&Purchasing probability&Price\\
     Portfolio optimization&Market environment&Asset returns&Allocation vector\\
     \bottomrule
    \end{tabular} 
\end{center}

\paragraph{Distributionally Robust Optimization (DRO)} 
DRO is a popular framework for solving optimization problems under uncertainty
that improves two classic approaches: 
Sample Average Approximation (SAA) \cite{shapiro2014lectures} and robust optimization \cite{bertsimas2004price, ben2009robust}.
SAA uses historical samples to construct and solve the ``empirical'' version of the objective function. 
Robust optimization seeks to optimize the worst-case scenario when the uncertain parameter is drawn from an uncertainty set.
The recent advances in DRO allow the decision-maker to incorporate richer information of the distribution into the optimization problem than both methods.
In particular, 
DRO aims to optimize the worst-case scenario when the \emph{distribution} of the uncertain parameter is drawn from an ambiguity set $\Fscr$, i.e. a family of distributions satisfying partial information inferred from the historical observations.
More precisely,
DRO can be reformulated as follows:
\begin{equation}\label{eq:dro}
\inf_{z \in \Zscr}\sup_{\mu_Y \in \mathcal{F}}\E_{\mu_Y}[c(z, Y)].
\end{equation}
The ambiguity set can be constructed based on the moments of the distribution or statistical distances to the empirical distribution. 
Both can be obtained in a data-driven manner.
The latter approach has gained popularity recently, including the Wasserstein metric \cite{esfahani2018data,gao2016distributionally}, Prohorov metric \cite{erdougan2006ambiguous}, goodness-of-fit test \cite{bertsimas2018robust} and  $\phi$-divergence \cite{ben2013robust}. We refer readers to a more comprehensive discussion in \cite{rahimian2019distributionally}.



We use DRO with the Wasserstein ambiguity set to solve \eqref{eq:formulation}, 
The statistical guarantees of the Wasserstein ambiguity set have been studied in  \cite{gao2016distributionally,gao2017distributional,blanchet2019quantifying}.
The computational tractability and the relationship to convex optimization also make it appealing \cite{zhao2014data, esfahani2018data,blanchet2019quantifying}. 
The Wasserstein metric has also been applied to statistical learning \cite{shafieezadeh2015distributionally} and operations management \cite{wozabal2014robustifying}.
However, the formulation \eqref{eq:dro} doesn't incorporate the covariate information.
DRO associated problems have been addressed in statistics and machine learning community in \cite{nguyen2020distributionallyEst, nguyen2020distributionallyMLE} and real applications into portfolio optimization \cite{nguyen2021robustifying}, single-stage stochastic programming \cite{dou2019distributionally} and quantile regression of fixed design \cite{qi2021distributionally} in recent years.
We consider the general cost function associated with prescriptive analytics. Compared with those application problems, this study is one of the first to investigate the theoretical properties and tractability when applying the Wasserstein metric to the DRO associated with prescriptive analytics \eqref{eq:formulation}.

Our framework is closely related to \cite{bertsimas2017bootstrap,kannan2020residuals}. 
\cite{kannan2020residuals} construct the Wasserstein ambiguity set around the residuals, while ours is on the distribution of $Y|X=x$. 
We provide a theoretical guarantee for the convergence of kernel estimation distribution instead of the bootstrap version in \cite{bertsimas2017bootstrap}.


\paragraph{Kernel Estimates}
To adapt the DRO framework and incorporate the covariate, we use kernel estimates, which is a classic approach in nonparametric statistics \cite{gyorfi2006distribution}.
More precisely, instead of centering the ambiguity set at the empirical distribution of $\{y_i\}_{i=1}^n$,
we use the historical samples of the covariate and construct a weighted empirical distribution of $\{y_i\}_{i=1}^n$, whose weight is proportional to the Nadaraya-Watson kernel estimator.
The intuition is clear: we weigh the historical samples more if the covariate is closer to the current covariate $x$.
The statistical properties of the kernel estimates have been studied thoroughly \cite{nadaraya1964estimating,watson1964smooth, tsybakov2008introduction}.
It has been applied to various areas in regression tasks due to its weak assumptions and flexibility for theoretical analysis. 
It has also been studied in different areas later such as stochastic optimization \cite{hannah2010nonparametric}, dynamic programming \cite{hanasusanto2013robust} and newsvendor problem \cite{ban2019big}. Various properties of asymptotic convergence and finite-sample guarantees have been attained in \cite{bertsimas2020predictive} with general assumptions under prescriptive analytics with real-data performance. To mitigate uncertainty in the stochastic setting, robust techniques have been studied in connection with the Nadaraya-Watson estimator for prescriptive analytics recently in \cite{ho2019data}, which incorporates the standard deviation term of stochastic objective as a regularization term to modify the original problem. 
Our work is different from the literature as 
we study the theoretical properties when applying the kernel estimates to the DRO setting.



\paragraph{Our Contributions}
In this paper, we propose a distributionally robust perspective on prescriptive analytics.
Motivated by the fact that robust techniques may improve the out-of-sample performance when the sample size is small, we adopt the nonparametric learning approach and construct Wasserstein-based ambiguity set.
To the best of our knowledge, no statistical guarantees have been established when the center of ambiguity set is not the empirical distribution with uniform weights among historical samples. 
In this paper, we aim to close this gap by incorporating the well-known Nadaraya-Watson estimator. 
We provide the measure concentration results for the nominal distribution in a wide range of Nadaraya-Watson kernels. 
Following that, we derive statistical properties including finite-sample performance guarantees and asymptotic convergences. 
Based on the ambiguity set, we show that the DRO can still be reformulated as convex problems to have computational tractability.
We summarize our main contributions as follows:

$\bullet$ \textbf{Distributionally robust prescriptive analytics using kernel estimates}. We provide a general formulation for a wide class of real-world problems using Nadaraya-Watson kernel estimators and DRO. 
Our data-driven approach is nonparametric, flexible to implement and requires very mild technical assumptions.

$\bullet~$\textbf{Measure Concentration and Statistical Guarantees}. For general Nadaraya-Watson kernels, we develop measure concentration results for the center of the ambiguity set and the actual unknown distribution, when the sample size grows. 
Because of the non-uniform weights, we extend the analysis in \cite{bolley2007quantitative, fournier2015rate}.
Based on the measure concentration results, we establish out-of-sample performance guarantees and asymptotic convergence results with regards to the underlying actual distribution for the DRO problem.

$\bullet$ \textbf{Problem Reformulation and Interpretation}. Leveraging techniques in Wasserstein-DRO reformulation, we show that strong duality holds in our problem.
As a result, it can be reformulated as convex problem. Furthermore, the problem is equivalent to the stochastic programming incorporating the weighted gradient regularization terms asymptotically. 
As two applications, we develop equivalent tractable reformulations for the newsvendor problem and CVaR-based portfolio optimization.
\paragraph{Technical notations} Throughout the paper, we use upper case letters to represent random quantities and lower case letters to represent deterministic quantities or samples.
We denote $[n] \coloneqq \{1,2,...,n\}$, $\|\cdot\|_p$ as the $\ell_p$-norm. 
We use $\delta_{x}$ to represent Dirac measure on the single point $  x$ and $\I{\cdot}$ denote the indicator function. 
Denote $\sigma_{\mathcal{Y}}(\cdot)$ as the support function of the domain $\mathcal{Y}$, i.e., $\sigma_\Yscr (\xi)=\sup_{y \in \Yscr}\{\xi^{T} y\}$. 
Denote $f^*(\cdot)$ as the convex conjugate of $f(\cdot)$ and $\|z\|_*$ denote the dual norm of $z$.
For a set $\Xscr$, denote $\mu(\Xscr)$ as the family of probability distributions on $\Xscr$.
We use $\mu_{Y|x}$ to denote the conditional probability measure of $Y$ given $X=x$; $\hat\mu_Y$ for a generic measure of $Y$ based on the historical samples, and specified to the kernel version $\mu_{n,x}$ in \eqref{eq:kernel}.
For the random sequences $\alpha_n, \beta_n$, we write $\alpha_n = o_p(\beta_n)$ if $\alpha_n/\beta_n$ converges to zero in probability or $\alpha_n = O_p(\beta_n)$ if $\alpha_n/\beta_n$ is uniformly bounded in probability. 
We use $\mathcal N(\mu,\sigma^2)$ as the normal distribution and $U(\cdot)$ as the uniform distribution.

\section{Model Formulation}\label{sec:formulation}
In this paper, we use the Wasserstein metric for the DRO framework.
\begin{definition}
\normalfont
(Wasserstein distance). The type-$p$ $(1 \leq p < + \infty)$ Wasserstein distance between two distributions $\mu$ and $\nu$ supported on $\Xi$ is defined as:
\begin{equation*}
        \Wscr_p(\mu,\nu) = \inf_{\xi \in \mu{(\Xi^2)}}\left\{\left(\int_{\Xi^2} \|x - y\|^p\xi(dx,dy)\right)^{\frac{1}{p}}: \xi_x=\mu, \xi_y=\nu\right\},
\end{equation*}
where $\xi_x$ and $\xi_y$ are the marginal distribution of $\xi$.
\end{definition} 
We first construct an ambiguity set for $\mu_Y$ based on the historical samples $\left\{(  x_i,y_i)\right\}_{i=1}^n$ and $X=x$.
If such a set is constructed, then we can apply the DRO framework \eqref{eq:dro} and solve 
\begin{equation}\label{eq:target}
    \hat{J}_n \coloneqq \inf_{z \in \mathcal{Z}}\sup_{\mu_Y\in \Fscr_n(x)}\E_{\mu_{Y}}[c(z,Y) | X = x].
\end{equation}

Specifically, we 
use the following ambiguity set in \eqref{eq:target} with the center $\hat{\mu}_Y$:
\begin{equation}\label{eq:ambiguity-set}
    \Fscr_n(x) = B_{\varepsilon}^p(\hat \mu_Y) \triangleq \left\{\mu: \Wscr_p\left(\hat\mu_Y,\mu\right)\leq\varepsilon\right\}.
\end{equation}
In other words, the ambiguity set is a ball centered at $\hat \mu_Y$, which depends on $X=x$ and the samples 
$\left\{(  x_i,y_i)\right\}_{i=1}^n$, with radius $\varepsilon $ defined using the Wasserstein distance.

Without the covariate information, the center $\hat\mu_Y$ is typically chosen as the empirical distribution of $\{Y_i\}_{i=1}^n$, i.e., $ \frac{1}{n}\sum_{i = 1}^{n}\I{Y = y_i}$.
Such problems are studied in \cite{esfahani2018data,gao2016distributionally}.
Failing to capture the extra information provided by the covariate, such a method does not lead to optimal decisions in our setting, even when the sample size $n$ tends to infinity.
For example, consider the objective in the newsvendor problem
\begin{equation*}
    c(z,Y) = b(Y - z)^+ + h(z-Y)^+,
\end{equation*}
for order quantity $z \in \R^+$, realized demand $Y$, unit backordering cost $b$ and unit holding cost $h$. 
If the covariate information is correlated with the demand 
\begin{equation*}
 Y| {X=x} \sim \beta^T   x + \varepsilon
\end{equation*}
for some mean-zero noise $\varepsilon$, then one can easily show that unless 
\begin{equation*}
\beta^T(  x - \mathbb{E}[X]) = 0,
\end{equation*}
ignoring the covariate doesn't lead to the optimal order quantity as $n\to\infty$.


Motivated by the Nadaraya-Watson estimator, we incorporate the historical samples and the current covariate 
by considering the following $\hat \mu_{Y}$:
\begin{equation}\label{eq:kernel}
    \hat\mu_Y = \mu_{n,x}(Y) \coloneqq \sum_{i=1}^{n} \frac{K( \frac{  x-  x_i}{h_n})}{\sum_{j=1}^n K( \frac{  x-  x_j}{h_n})} \I{Y = y_i}.
\end{equation}
The kernel function $K:\R^{d_x} \to [0,+\infty)$ and $h_n\in \R^+$ is the bandwidth that depends on the sample size $n$.
Usually $K(x)$ is large when $\|x\|_2$ is small, creating a local averaging estimate.
Typical choices of $K(x)$ include $\I{\|x\|_2\le 1}$, $(1-x^2)_+$ and $\exp(-x^2)$.
By using \eqref{eq:kernel} as the center of the Wasserstein ball, 
we capture the following pattern of the historical data:
when a sample $(x_i,y_i)$ has a covariate $x_i$ close to $x$, the kernel function $K((x-x_i)/h_n)$ tends to be large and the distribution $\hat \mu_Y$ assigns more weight to $y_i$.
It is a nonparametric way to capture the join distribution of $(X, Y)$ in the historical data and leverage the information $X=x$.
Note that solving \eqref{eq:target} using \eqref{eq:ambiguity-set} and \eqref{eq:kernel}
reduces to the problem studied in \cite{bertsimas2020predictive} if we set $\varepsilon = 0$.

\section{Theoretical Analysis}\label{sec:theoretical}
In this section, we provide the theoretical analysis of the formulation in Section~\ref{sec:formulation}.
We first give a few standard assumptions for the analysis.
\begin{assumption}\label{asp:generate}
\normalfont
The samples $\{(x_i,y_i)\}_{i=1}^n$ are generated independently from the following distribution:
$X_i\sim \mu_X$ in the domain $[0,1]^{d_x}$ and $Y_i=f(X_i)+\varepsilon_i$, where $\varepsilon$ are i.i.d. mean-zero noise and $\E[\varepsilon^2] \leq \sigma^2$.
Furthermore, $f$ is Lipschitz continuous, i.e. there exists $0 \le L < \infty$ such that $|f(  x)-f(  x_0)|\le L\|  x-  x_0\|_2$.
\end{assumption}
We denote the measure $Y$ given $X=x$, which is $f(x)+\varepsilon$, as $ \mu_{Y | x}$.
It is simply a shift of $\mu_\varepsilon$. 
Note that we only consider $d_y=1$ for simplicity, as the conclusion still holds for $d_y>1$.
Assumption~\ref{asp:generate} specifies the generative model of the historical data, which is standard in the literature.

\begin{assumption}\label{asp:distribution}
\normalfont
    We have $\mu_x(\cdot)\ge c>0$ and $\E[|\varepsilon|^{2p+\delta}]<\infty$ for some $c,\delta>0$ for type-$p$ Wasserstein distance.
\end{assumption}
The lower bound of the measure guarantees that for any covariate $x$, there are enough samples around $x$ when the sample size tends to infinity. 
The second part is a necessary condition to bound the Wasserstein distance, similar to \cite{bolley2007quantitative,fournier2015rate}.

The next assumption is concerned with the kernel function for the Nadaraya-Watson estimator.
\begin{assumption}\label{asp:kernel}
\normalfont
Kernel function $K(\cdot)$ is decreasing in $[0, \infty)$ and there exist $R\geq r>0$ and $b_R\geq b_r>0$ such that $K(\cdot)$ satisfies
\begin{equation*}
    b_r \I{\|  x\|_2\le r}\le K(  x)\le b_R \I{\|  x\|_2\le R}.
\end{equation*}
\end{assumption}
In other words, the kernel function is upper- and lower-bounded by two naive kernels.
It is a technical assumption to guarantee the statistical properties. 

\subsection{Measure Concentration}
We first establish the measure concentration under the Wasserstein distance. 
This serves as a building block for the subsequent statistical analysis.
\begin{theorem}\label{prop:ws-concentration}
\normalfont
Suppose the kernel function in \eqref{eq:kernel} satisfies Assumption \ref{asp:kernel}.
For all $ x\in [0,1]^{d_x}$, $t\in [0,1]$ and a sequence of $\{h_n\}_{n=1}^{+\infty}$ such that $h_n\to0$ and $nh_n^{d_x}\to\infty$,
we have
\begin{equation}\label{eq:convergence_rate}
    \mathbb{P}(\Wscr_p( \mu_{Y | x},\mu_{n,x})>t)\le c_1 \exp(-c_2n t^2),
\end{equation}
where $c_1$ and  $c_2$ are independent of $n$, $t$ and $  x$.
\end{theorem}

To prove the theorem, we first investigate the case when $K(x)$ represents the naive kernel (also referred to as the window/box kernel), i.e. $K(x) = \I{\|x\|_2\le r}$.
That is,
\begin{equation*}
\mu^r_{n,x}\coloneqq \sum_{i=1}^n \frac{ \I{\|  x-  x_i\|_2\le rh_n}}{\sum_{j=1}^n \I{\|  x-  x_j\|_2\le rh_n}} \I{Y = y_i}.
\end{equation*}
Because of Assumption~\ref{asp:generate}, we can consider the shifted measure
\begin{equation}\label{eq:approximation_measure}
\bar\mu^r_{n,x}:= \sum_{i=1}^n \frac{ \I{\|  x-  x_i\|_2\le rh_n}}{\sum_{j=1}^n \I{\|  x-  x_j\|_2\le rh_n}} \I{Y = y_i-f(  x_i)+f(  x)}
\end{equation}
The benefit of analyzing $\bar\mu_{n,x}^r$ is that after the shift in $y_i$, $y_i-f(x_i)+f(x)$ can be regarded as i.i.d. samples.
This allows us to leverage the results in \cite{fournier2015rate} to control its Wasserstein distance from $\mu_{Y|x}$, by focusing on the samples in the ball of radius $rh_n$.
We can then control the Wasserstein distance between $\bar\mu^r_{n,x}$ and $\mu_{n,x}^r$ using the Lipschitz continuity of $f(x)$.

Next, using the result for naive kernels and Assumption~\ref{asp:kernel}, we argue that 
\begin{equation*}
    \Wscr( \mu_{Y | x}, \mu_{n, x}) \leq \max\left\{ \Wscr( \mu_{Y | x},\mu_{n,x}^+), \Wscr( \mu_{Y | x},\mu_{n,x}^-)\right\},
\end{equation*}
where $\mu_{n,x}^+$ and $\mu_{n,x}^-$ are the measures associated with the two kernel bounds in Assumption~\ref{asp:kernel}.
The complete proof can be found in the appendix.
\begin{remark}
\normalfont
When $d_y>1$, the RHS of \eqref{eq:convergence_rate} changes to $c_1 \exp(-c_2 n t^{d_y})$. 
This is the situation for the numerical experiments in Section~\ref{sec:portfolio}.
\end{remark}

Theorem~\ref{prop:ws-concentration} provides guidelines to the choice of $\varepsilon$ in \eqref{eq:ambiguity-set}, which is the radius of the Wasserstein ball. 
In particular, 
for $\alpha \in(0,1)$, if $n \geq \frac{\log(\frac{c_1}{\alpha})}{c_2}$, then a Wasserstein ball centered at $\mu_{n,x}$ with radius $\varepsilon_n(\alpha)\coloneqq \sqrt{\frac{\log(\frac{c_1}{\alpha})}{c_2n}}$ contains the actual measure $\mu_{Y|x}$ with probability at least $1-\alpha$. 
When we fix $\alpha$, the radius of the ambiguity set $\varepsilon_n(\alpha) \to 0$ as $n \to \infty$. 
Let $\hat z_n$ be the optimal solution to \eqref{eq:target}.
The next result states that the DRO formulation leads to a solution that performs well under the actual distribution when compared to $\hat J_n$.
This is referred to as the out-of-sample guarantee for the objective function.
\begin{proposition}\label{prop:sampleguarantee}
\normalfont
    Suppose Assumptions~\ref{asp:generate}, \ref{asp:distribution} and~\ref{asp:kernel} hold. 
    We have
    \begin{equation}
        \mathbb P\left(\E_{ \mu_{Y | x}}[c(\hat{z}_n,Y)]\leq \hat{J}_n\right) \geq 1-\alpha.
    \end{equation}
\end{proposition}

\subsection{Convergence of the Objective Function}
We next study the convergence of the optimal value of \eqref{eq:target} to the actual optimal value given $X=x$:
\begin{equation}\label{eq:true_target}
    J^*\coloneqq \inf_{z \in \mathcal{Z}}\mathbb{E}_{ \mu_{Y | x}}[c(z, Y)]. 
\end{equation}
We would impose the following regular assumption following \cite{esfahani2018data}.
\begin{assumption}\label{asp:cost_function}
\normalfont
Given $z$, $c(z, y)$ is proper, Lipschitz continuous with constant $L_z < \infty$ and concave for $y\in \Yscr$, where $\Yscr$ in a closed convex set. 
Given $y$, $c(z, y)$ is upper semi-continous in $z$.
Furthermore, there exists $0 \leq L_{c, p} < \infty$ such that
$|c(z,y)| \leq L_{c,p}(1 + \|y\|^p), \forall z \in \mathcal{Z}, y \in \mathcal{Y}$.
\end{assumption}
This assumption holds naturally for most functions used in real applications. 
It allows us to derive the following result.
\begin{proposition}\label{prop:asymptotic_convergence}
\normalfont
Suppose Assumptions~\ref{asp:generate}, \ref{asp:distribution}, \ref{asp:kernel} and~\ref{asp:cost_function} hold.
If the sequence $\{\alpha_n\}_{n = 1}^{+\infty}$ satisfies $\alpha_n \in (0,1)$, $\sum_{n = 1}^{+\infty}\alpha_n < \infty$, and $\lim_{n \to \infty}\varepsilon_n(\alpha_n) = 0$, then we have $\lim\inf_{n\to \infty}\hat{J}_n= J^*$ almost surely.
In addition, if $\mathcal{Z}$ is closed, then any accumulation point of $\{\hat{z}_n\}_{n = 1}^{+\infty}$ for \eqref{eq:target} converges to an optimal solution $z^*$ to \eqref{eq:true_target} almost surely.
\end{proposition}

Proposition~\ref{prop:asymptotic_convergence} guarantees that when the historical sample size increases, 
the formulation \eqref{eq:target} converges to \eqref{eq:true_target} in terms of the optimal solution and optimal value.
We can further characterize the convergence rate below.
\begin{proposition}\label{prop:convergence_bound}
\normalfont
Suppose Assumptions~\ref{asp:generate}, \ref{asp:distribution}, \ref{asp:kernel} and~\ref{asp:cost_function} hold.
If  $c(z, \cdot)$ is Lipschitz continuous in $z$, then we have
\begin{equation}
    |\hat{J}_n - J^*| = O_p(\varepsilon_n(\alpha_n)), \quad |\mathbb{E}_{\mu_{Y|x}}[c(\hat{z}_n, Y)] - J^*| = O_p(\varepsilon_n(\alpha_n)).
\end{equation}
\end{proposition}
Therefore, our formulation using the kernel functions allows to reduce the DRO problem with covariates to the standard one studied in 
\cite{esteban2020distributionally} and most statistical properties still hold. 

\section{Problem Reformulation}\label{sec:reformulation}
In this section, we provide a reformulation of \eqref{eq:target} for the computational tractability.
\subsection{Duality}
We first define the primal and dual problems:
\begin{align}
    c_P(z)& \coloneqq \sup_{\mu_Y\in B_{\varepsilon}^p(\hat\mu_{Y})}\E_{\mu_{Y}}\left[c(z,Y)\right],\notag\\
    c_D(z)& \coloneqq \inf_{\lambda \geq 0}\{\lambda \varepsilon^p - \E_{\xi \sim {\hat{\mu}_Y}}\left\{\inf_{y \in \mathcal{Y}}[\lambda \|\xi - y\|_p^p- c(z, y)]\right\}, \forall z\in \mathcal{Z}. 
\label{eq:dual}
\end{align}
we will show $c_P(z) \leq c_D(z)$ in the appendix along with the strong duality property later.
Using \eqref{eq:kernel} as the nominal distribution with the Nadaraya-Watson estimator, the primal and dual problems can be simplified to
\begin{equation*}
c_P(z) = c_D(z)= \min_{\lambda \geq 0}\left\{\lambda \varepsilon^p + \sum_{i = 1}^{n}w_i\sup_{y \in \mathcal{Y}}[c(z,y) - \lambda \|y- y_i\|_p^p]\right\},    
\end{equation*}
where $w_i = \frac{K( \frac{ x-  x_i}{h_n})}{\sum_{j=1}^n K( \frac{  x-  x_j}{h_n})}$
and the dual problem \eqref{eq:dual} has a minimizer $\lambda^*$.
Strong duality can then be established.
\begin{lemma}\label{lemma:duality}
\normalfont
Suppose Assumption \ref{asp:cost_function} holds. 
We have
\begin{equation}\label{eq:duality2}
c_P(z) = c_D(z) =\sup_{\xi_i \in \mathcal{Y}}\left\{\sum_{i = 1}^{n}w_i c(z, \xi_i): \sum_{i = 1}^{n}w_i\|\xi_i- y_i\|_p^p \leq \varepsilon^p \right\}.
\end{equation}
Simultaneously, the worst-case distribution is supported on at most $n + 1$ points and has the form
\begin{equation*}
    \tilde{\mu}_Y \coloneqq \sum_{i \not=i_0}w_i\delta_{y^*_i} + w^*\delta_{y^*_{1i_0}} + (w_{i_0}-w^*)\delta_{y^*_{2i_0}},
\end{equation*}
where $i_0\in \{1,\dots,n\}$ and $w^*\in [0, w_{i_0}]$.
The support $y_i^*$, $i\neq i_0$, satisfies 
$y^*_i \in \argmin_{y \in \mathcal{Y}}\{\lambda^*\|y- y_i\|_p^p- c(z,y)\}$; 
similarly, $y^*_{1i_0}, y^*_{2i_0} \in \argmin_{y \in \mathcal{Y}}\{\lambda^*\|y- y_{i_0}\|_p^p- c(z,y)\}$.
\end{lemma}

\subsection{Regularization Interpretation}
We now provide an alternative perspective for our formulation using regularization.
\begin{proposition}\label{prop:regularization}
\normalfont
For any $z$, when $c(z,y)$ satisfies assumptions in Theorem 2 of \cite{gao2017distributional}, the following condition holds almost surely:
$$\left|\sup_{\mu_Y \in B_{\varepsilon_n(\alpha_n)}^p(\hat\mu_Y)} \mathbb{E}_{\mu_Y}[c(z, Y)] - \left(\mathbb{E}_{\hat{\mu}_Y}[c(z, Y)] + \alpha_n\|\nabla_y c(z, y)\|_{\hat{\mu}_Y, p^*}\right)\right| = o_p(\alpha_n), $$
where the gradient-norm regularization is defined as follows:
\begin{equation*}
    \|\nabla_y c(z, y)\|_{\hat{\mu}_Y, p}:= 
    \begin{cases}
    (\sum_{i \in [n]}w_i\|\nabla_y c(z, y_i)\|_p^{p})^{\frac{1}{p}},&1\leq p < +\infty, \\
    \max_{i \in [n]}\|\nabla_y c(z, y_i)\|_p,&p = \infty.
    \end{cases}
\end{equation*}
\end{proposition}
In other words, the formulation \eqref{eq:target} with the Wasserstein ball centered at \eqref{eq:kernel} is penalizing the weighted gradients at the historical samples.
Compared with the DRO problem without covariates in \cite{gao2017distributional}, 
the weights $w_i$ on the history samples capture the sample ``importance'' depending on the distance from the current covariate given by the kernel functions.

\subsection{Computational Tractability}
We consider $p=1$ in the Wasserstein distance in this section for the ambiguity set.
For the numerical studies in Section~\ref{sec:numerical}, the objective functions may be piece-wise linear such as the newsvendor problem.
To accommodate such a situation, we impose an assumption below.
\begin{assumption}\label{asp:cost_function_piecewise}
\normalfont
Given $z$, $c(z,y)$ can be written as $ \max_{k \in [K]}c_k(z, y)$, abbreviated as $c_{z, k}(y)$, where each constituent function $c_k(z, y)$ is proper, concave and upper semi-continuous in $y$.
\end{assumption}
\begin{proposition}\label{prop:reformulation}
\normalfont
Suppose Assumption \ref{asp:cost_function_piecewise} holds. For any $\varepsilon>0$ in the Wasserstein ambiguity set, we have
\begin{equation}\label{eq:reformulation}
\begin{aligned}
c_P(z) = \inf_{\lambda, s_i, p_{ik}, v_{ik}}~&\lambda \varepsilon + \sum_{i = 1}^{n}w_is_i\\
\text{s.t.}~&[-c_{z, k}]^*(p_{ik}- v_{ik}) + \sigma_{\mathcal{Y}}(v_{ik}) - p_{ik}y_i \le s_i, \forall i \in [n], \forall k \in [K]\\
~&|p_{ik}| \le \lambda, \forall i \in [n], k \in [K].\\
\end{aligned}
\end{equation}
\end{proposition}
We refer to the introduction for the notations in the formulation.
In our numerical studies, Proposition~\ref{prop:reformulation} ensures that the problems can be reformulated into convex optimization problems and guarantees computationally tractability.

\section{Numerical Studies}\label{sec:numerical}
Next, we demonstrate the performance of our framework in synthetic and empirical experiments. 
All optimization problems are implemented in Python via the convex optimization solver Gurobi with an Intel(R) Core(TM) i7-8650U CPU @1.90GHz Personal Computer.

\paragraph{Benchmarks} We denote our model as \texttt{NW-DRO} and optimize via Gaussian kernels\footnote{We also test naive kernels and results are almost the same.}. 
Additionally, we compare the performance with several benchmarks: \texttt{Naive-SO} represents the SAA method; \texttt{Naive-DRO} represents the distributionally robust optimization without covariates \cite{esfahani2018data}; \texttt{NW-SO} represents the prescriptive analytics under Nadaraya-Watson estimator in \cite{bertsimas2020predictive}, 
which is equivalent to $\epsilon=0$ for the ambiguity set.

\subsection{Newsvendor Problem}
The objective of the newsvendor problem  $c(z, Y) = \max{b(Y - z), h(z - Y)}$, which is piece-wise linear.
We use $b = 10$ and $h = 1$ across our experiments. 
The historical data is generated similarly to the setting in \cite{bertsimas2017bootstrap}.
The covariate includes the temperature $t$ ($^{\circ}C$) and the day in a week $d$, generated from the distribution $t \sim \mathcal{N}(20, 4)$ and $d\sim U(1,2,\dots,7)$.
The conditional distribution of the demand $Y$ given $x=(t,d)$ is $\mathcal{N}(100 + (t - 20) + 20 \cdot \I{d \in \{\text{Weekend}\}}, 16)$.

We focus on the measure ``disappointment'', defined as $\mathbb P(\mathbb{E}_{\mu_{Y|x}}[c(\hat{z}_n, Y)] \geq \hat{J}_n)$, by averaging 2500 simulated instances in the out-of-sample tests.
A lower disappointment indicates better performance.
The measure has also been used in \cite{bertsimas2017bootstrap, esfahani2018data, van2020data}.
\begin{figure}[h]  
    \centering    
    \subfloat[Varying size with $C = 40$] 
    {
        \begin{minipage}[t]{0.5\textwidth}
            \centering          
            \includegraphics[width=0.9\textwidth, trim = 35 35 45 45]{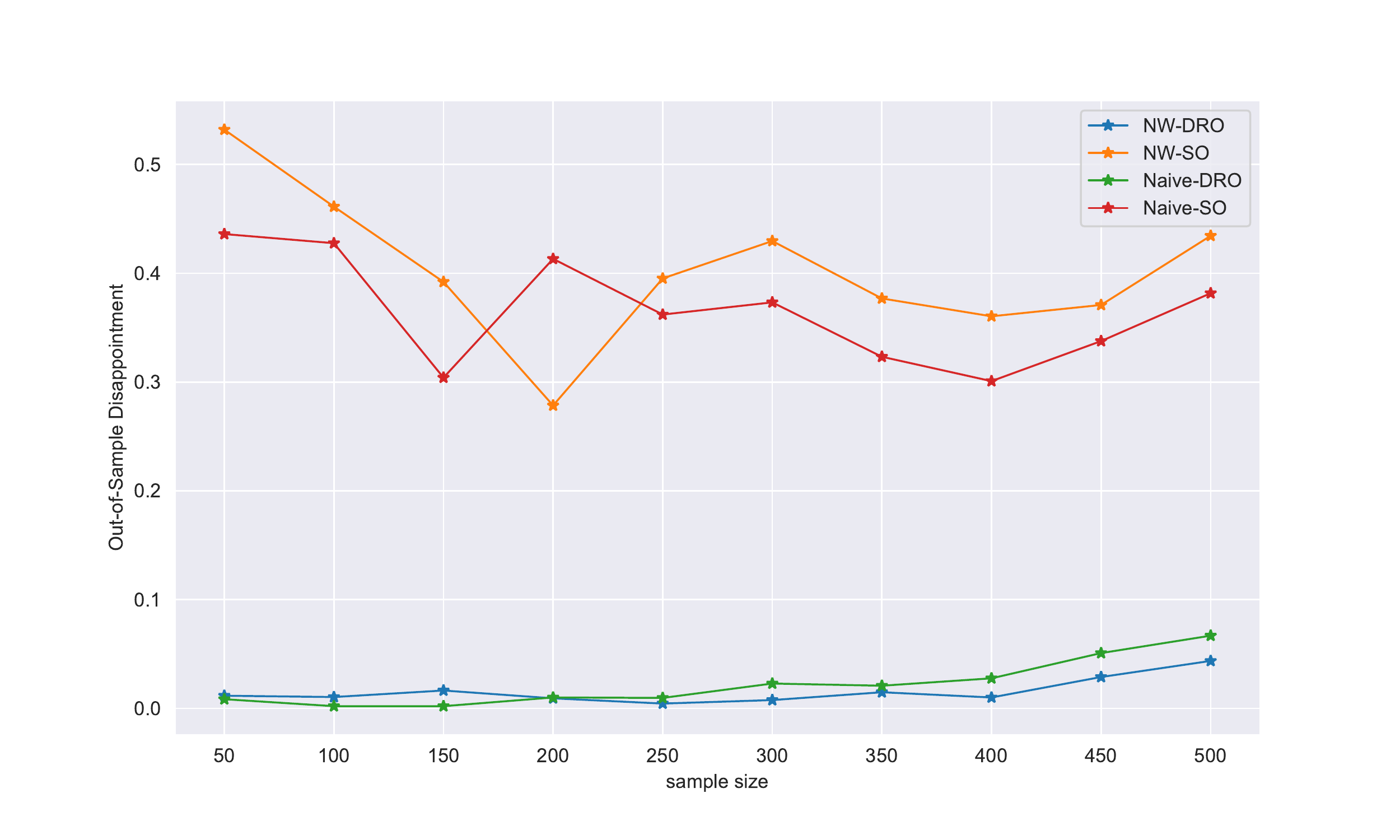}   
        \end{minipage}%
    }
    \subfloat[Varying size with different $C$] 
    {
        \begin{minipage}[t]{0.5\textwidth}
            \centering      
            \includegraphics[width=0.9\textwidth, trim = 35 35 45 45]{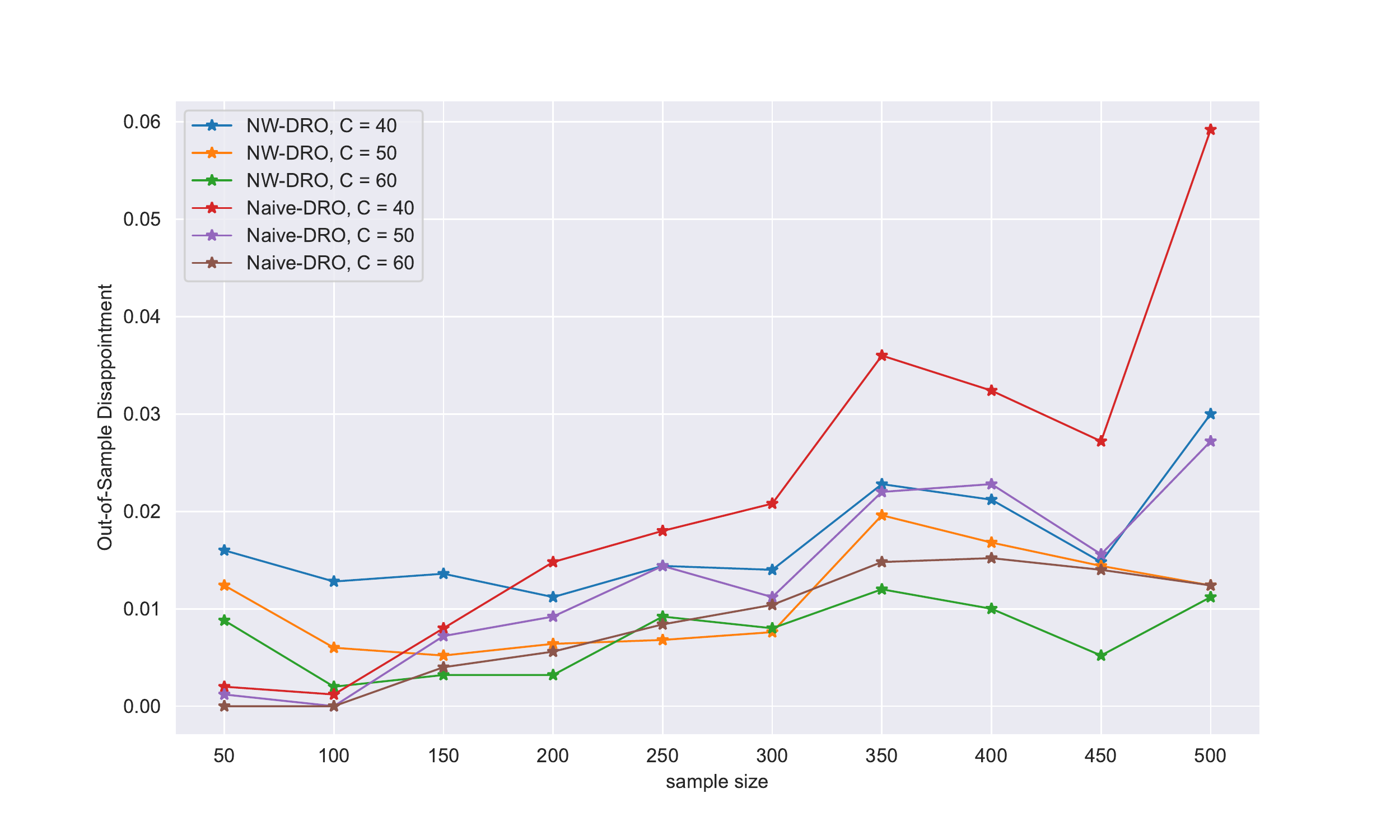}   
        \end{minipage}
    }%
    
    \subfloat[Fixed size with $\varepsilon$ = 2] 
    {
        \begin{minipage}[t]{0.5\textwidth}
            \centering          
            \includegraphics[width = 0.9\textwidth, trim = 35 35 45 45 ]{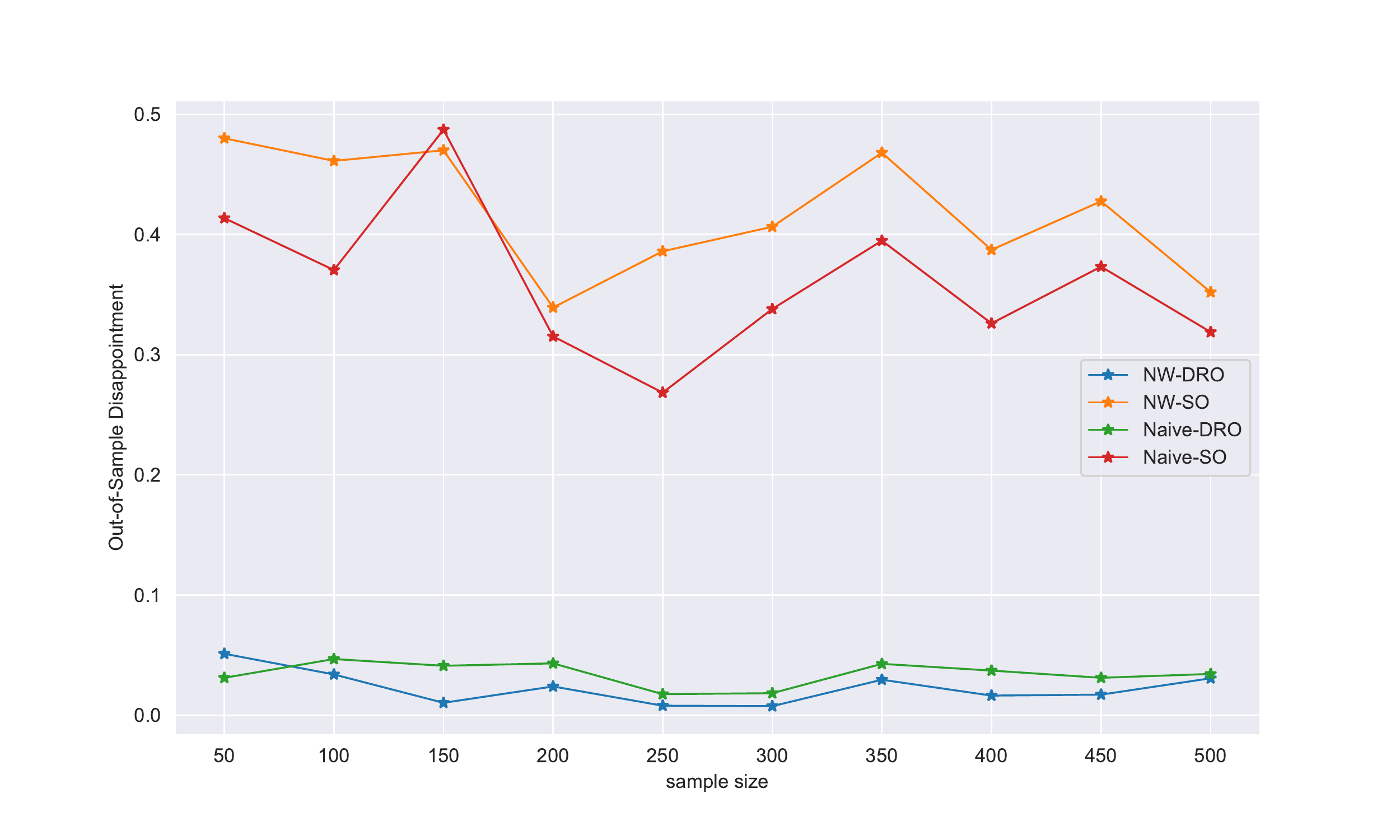}   
        \end{minipage}%
    }
    \subfloat[Fixed size with different $\varepsilon$] 
    {
        \begin{minipage}[t]{0.5\textwidth}
            \centering      
            \includegraphics[width = 0.9\textwidth, trim = 35 35 45 45]{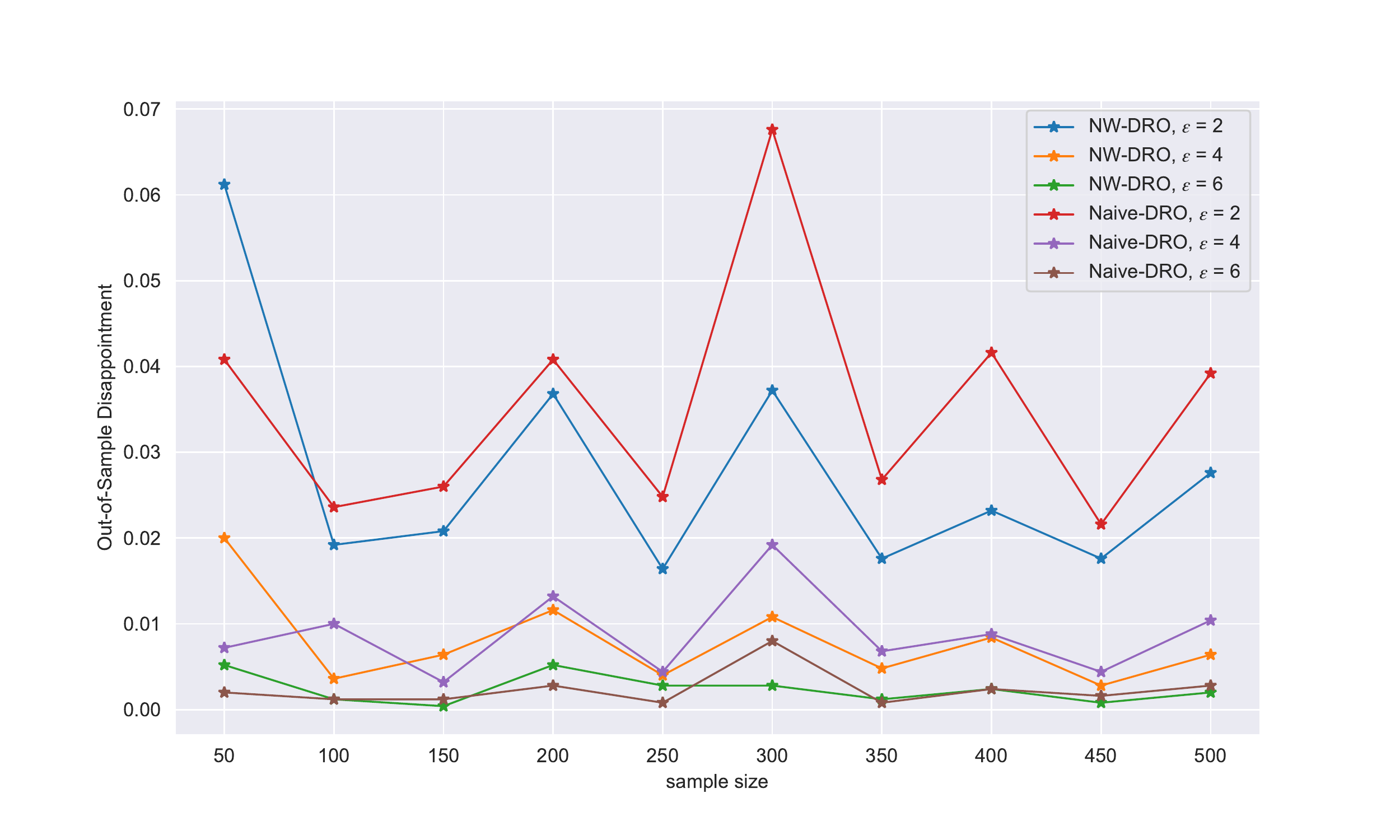}   
        \end{minipage}
    }%
    \caption{Out-of-sample disappointment against sample size $n$ based on 2500 simulations.
    The radius of the Wasserstein ambiguity set is specified by $\varepsilon$ (fixed size) or $\varepsilon=C/n$ (varying size). } 
    \label{fig1}  
\end{figure}

We investigate the performance when the ambiguity set has a fixed radius $\varepsilon$ or a sample-dependent radius $\varepsilon = \frac{C}{n}$ in Figure \ref{fig1}.
From the experiment, we can see that the DRO formulation significantly improves SO (\texttt{Naive-SO} versus \texttt{Naive-DRO} and \texttt{NW-SO} versus \texttt{NW-DRO}).
The formulation that incorporates the covariate only slightly improves the performance (\texttt{Naive-SO} versus \texttt{NW-SO} and \texttt{Naive-DRO} versus \texttt{NW-DRO}), which may be caused by the small sample size that is not sufficient to capture the complex joint distribution of $(X,Y)$.
Overall, our method performs the best.


\subsection{Portfolio Allocation}\label{sec:portfolio}
Using the real datasets about asset returns and covariates from the data library of Kenneth French's website\footnote{\url{http://mba.tuck.dartmouth.edu/pages/faculty/ken.french/data_library.html}},
we apply our framework to portfolio allocation.  
We use $  Y \in \mathbb{R}^{d_y}$ to represent the random asset returns and $z \in \mathcal{Z}= \{z \in \mathbb{R}^{d_z} | \sum_{i \in [d_z]} z_i = 1, z_i \ge 0\}$ for the allocation vector, where $d_y = d_x$. 
We focus on the following objective functions:
\begin{equation}\label{eq:cvarformulation}
\text{CVaR}_{\eta}(  R) = \inf_{v}\{v + \frac{1}{\eta}\mathbb{E}_{  R}[-   R - v]^+\}
\end{equation}
\begin{equation}\label{eq:cvar_objective}
c(z, Y) = \text{CVaR}_{\eta}(  Y^\top  z) - \gamma   Y^\top  z.
\end{equation}
The measure CVaR captures the trade-off between risk and the mean return.
In this experiment, we take $\eta = 0.05, \gamma= 1$. 

Using the linear reformulation of CVaR following \cite{rockafellar2002conditional}, 
it is easy to verify that the objective function satisfies Assumption~\ref{asp:cost_function_piecewise} by 
\begin{equation*}
\min_{z}\sup_{\mu_Y \in B_{\varepsilon}^p(\hat{\mu}_Y)}\E_{\mu_Y}[c(z, Y)] = \min_{z, v}\sup_{\mu_Y \in B_{\varepsilon}^p(\hat{\mu}_Y)}\E[\max\{-(\gamma + \frac{1}{\eta} Y^{\top}z) + (1 - \frac{1}{\eta}) v, -\gamma Y^{\top} z + v\}]    
\end{equation*}
Aside from the benchmarks mentioned above, we also investigate the Equally Weighted (\texttt{EW}) Model, i.e. set  $z_i = \frac{1}{d_z}$ for all $i \in [d_z]$. 
This model gains popularity in \cite{demiguel2007optimal} for its robustness. We take the most recently (i.e. last month) observed Fama-French three-factor model \cite{fama1992cross} of an associated portfolio monthly return data as the covariates.

To obtain and evaluate the out-of-sample performance, we apply the ``rolling-sample'' approach following \cite{demiguel2007optimal} on the monthly data from July 1963 to December 2018 ($T = 666$) with the estimation window to be $M = 60$ months\footnote{Rolling-sample: For data spanning $T$ months, given month $t$, we use the dataset from months $t$ to $t + M - 1$ as the input to solve the portfolio optimization problem. Then we use the optimal weight to compute the returns in month $t + M$. We repeat this procedure by adding the next and dropping the earliest month until $t = T - M$. This gives us $T - M$ monthly out-of-sample returns.}. We report the experimental results of different models with three measures: Sharpe Ratio (SR), Empirical CVaR (E-CVaR), and Certainty-equivalent Return (CEQ). E-CVaR is computed by \eqref{eq:cvarformulation} replacing $\mathbb{E}_{R}$ with the empirical average. SR and CEQ are computed with the formula in \cite{demiguel2007optimal}:
\begin{align*}
   \text{SR}&\coloneqq \widehat{\text{mean}}(\{  r_i\}_{i \in [T - M]})/\widehat{\text{std}}(\{  r_i\}_{i \in [T - M]}),\\
   \text{CEQ}&\coloneqq \widehat{\text{mean}}(\{  r_i\}_{i \in [T - M]}) - (\widehat{\text{std}}(\{  r_i\}_{i \in [T - M]}))^2.
\end{align*}
We set $\varepsilon_n = k M^{\frac{1}{d_y}}$ as the radius for the ambiguity set in the DRO setting, where $k$ is a hyperparameter to control.
The results are reported in Table~\ref{table:empirical}. 
\begin{table}[H]
    \centering
    \small
    \begin{tabular}{c|cccccc|ccc}
        \toprule
	    &\texttt{EW}&\texttt{Naive-SO}&\texttt{NW-SO}&\multicolumn{3}{c}{\texttt{Naive-DRO}}&\multicolumn{3}{|c}{\texttt{NW-DRO}}\\
		\hline
        $k$&-&-&-&0.2&0.4&0.8&0.2&0.4&0.8\\
        \hline
        \multicolumn{10}{c}{6 Portfolios Formed on Size and Book-to-Market (2 x 3)}\\
        \hline
        SR     & 0.2045 & 0.2070   & 0.2130 & 0.2045    & 0.2045 & 0.2045 &\textbf{0.2243}& 0.2217 & 0.2112  \\
        E-CVaR & 0.1089 &\textbf{0.0952}& 0.0954 & 0.1089    & 0.1089 & 0.1089 & 0.1008 & 0.1039 & 0.1075  \\
        CEQ    & 0.0076 & 0.0072   & 0.0075 & 0.0076    & 0.0076 & 0.0076 &\textbf{0.0082}& 0.0081 & 0.0079 \\
        \hline
        
        \multicolumn{10}{c}{10 Industry Portfolios}\\
        \hline
        SR     & 0.2203 & 0.2288 & 0.2301   & 0.2213    & 0.2203 & 0.2203 & 0.2535 &\textbf{0.2550}& 0.2385  \\
        E-CVaR & 0.0917 & 0.0812 & 0.0800   & 0.0909    & 0.0917 & 0.0917 &\textbf{0.0793}& 0.0829 & 0.0878  \\
        CEQ    & 0.0075 & 0.0072 & 0.0072   & 0.0075    & 0.0075 & 0.0075 & 0.0081 &\textbf{0.0083}& 0.0080 \\
        \hline
        \multicolumn{10}{c}{25 Portfolios Formed on Book-to-Market and Investment (5 x 5)}\\
        \hline
        SR     & 0.2029 & 0.2302   & 0.2265 & 0.2225    & 0.2119 & 0.2045 & 0.2292 & 0.2292 &\textbf{0.2321}\\
        E-CVaR & 0.1121 & 0.0962  &\textbf{0.0959}& 0.1068    & 0.1101 & 0.1112 & 0.0969 & 0.0984 & 0.1011  \\
        CEQ& 0.0078 & 0.0081   & 0.0081 & 0.0084    & 0.0081 & 0.0078 & 0.0082 & 0.0083 &\textbf{0.0086} \\
        \hline
        \multicolumn{10}{c}{30 Industry Portfolios}\\
        \hline
        SR     & 0.1981 & 0.2421   & 0.2350 & 0.2220    & 0.2067 & 0.2013 & 0.2301 & 0.2400 &\textbf{0.2430}\\
        E-CVaR & 0.1042 &\textbf{0.0846} & 0.0872 & 0.0927    & 0.0990 & 0.1033 & 0.0871 & 0.0850 & 0.0871  \\
        CEQ    & 0.0072 & 0.0079   & 0.0076 & 0.0077    & 0.0073 & 0.0072 & 0.0074 & 0.0079 &\textbf{0.0082} \\
		\hline
		\bottomrule
\end{tabular}
    \caption{Comparison of different policies in the portfolio allocation problem.}
    \label{table:empirical}
\end{table}
From Table \ref{table:empirical}, we can see that
our \texttt{NW-DRO} method outperforms other policies in terms of SR and CEQ. 
It performs well in terms of the empirical CVaR.
For higher dimension $d_y$, the optimal choice of the size of the ambiguity set tends to be larger to incorporate the higher estimation error.


\section{Conclusion and Future Research} \label{sec:conclusion}
We develop a distributionally robust framework that uses the Nadaraya-Watson kernel estimator as the nominal distribution in the Wasserstein ambiguity set to incorporate the covariate. 
We establish theoretical properties on the out-of-sample performance and asymptotic convergence and validate the practical performance in real applications.

It remains an open question to extend the statistical analysis beyond kernels with finite support (Assumption~\ref{asp:kernel}).
Another exciting direction is to develop a framework that can identify redundant variables when the covariate is sparse.
It would require a nonparametric variable selection procedure \cite{li2020dimension} that can be incorporated in the construction of the ambiguity set.
\bibliographystyle{unsrt}

\begin{thebibliography}{10}
	
	\bibitem{elmachtoub2017smart}
	Adam~N Elmachtoub and Paul Grigas.
	\newblock Smart ``predict, then optimize''.
	\newblock {\em arXiv preprint arXiv:1710.08005}, 2017.
	
	\bibitem{tulabandhula2013machine}
	Theja Tulabandhula and Cynthia Rudin.
	\newblock Machine learning with operational costs.
	\newblock {\em J Mach Learn Res}, 14(1):1989--2028, 2013.
	
	\bibitem{bertsimas2020predictive}
	Dimitris Bertsimas and Nathan Kallus.
	\newblock From predictive to prescriptive analytics.
	\newblock {\em Management Science}, 66(3):1025--1044, 2020.
	
	\bibitem{bertsimas2021data}
	Dimitris Bertsimas and Nihal Koduri.
	\newblock Data-driven optimization: A reproducing kernel hilbert space
	approach.
	\newblock {\em Operations Research}, 2021.
	
	\bibitem{notz2021prescriptive}
	Pascal~M Notz and Richard Pibernik.
	\newblock Prescriptive analytics for flexible capacity management.
	\newblock {\em Management Science}, 2021.
	
	\bibitem{shapiro2014lectures}
	Alexander Shapiro, Darinka Dentcheva, and Andrzej Ruszczy{\'n}ski.
	\newblock {\em Lectures on stochastic programming: modeling and theory}.
	\newblock SIAM, 2014.
	
	\bibitem{bertsimas2004price}
	Dimitris Bertsimas and Melvyn Sim.
	\newblock The price of robustness.
	\newblock {\em Operations Research}, 52(1):35--53, 2004.
	
	\bibitem{ben2009robust}
	Aharon Ben-Tal, Laurent El~Ghaoui, and Arkadi Nemirovski.
	\newblock {\em Robust optimization}.
	\newblock Princeton University Press, 2009.
	
	\bibitem{esfahani2018data}
	Peyman~Mohajerin Esfahani and Daniel Kuhn.
	\newblock Data-driven distributionally robust optimization using the
	wasserstein metric: Performance guarantees and tractable reformulations.
	\newblock {\em Mathematical Programming}, 171(1):115--166, 2018.
	
	\bibitem{gao2016distributionally}
	Rui Gao and Anton~J Kleywegt.
	\newblock Distributionally robust stochastic optimization with wasserstein
	distance.
	\newblock {\em arXiv preprint arXiv:1604.02199}, 2016.
	
	\bibitem{erdougan2006ambiguous}
	Emre Erdo{\u{g}}an and Garud Iyengar.
	\newblock Ambiguous chance constrained problems and robust optimization.
	\newblock {\em Mathematical Programming}, 107(1):37--61, 2006.
	
	\bibitem{bertsimas2018robust}
	Dimitris Bertsimas, Vishal Gupta, and Nathan Kallus.
	\newblock Robust sample average approximation.
	\newblock {\em Mathematical Programming}, 171(1):217--282, 2018.
	
	\bibitem{ben2013robust}
	Aharon Ben-Tal, Dick Den~Hertog, Anja De~Waegenaere, Bertrand Melenberg, and
	Gijs Rennen.
	\newblock Robust solutions of optimization problems affected by uncertain
	probabilities.
	\newblock {\em Management Science}, 59(2):341--357, 2013.
	
	\bibitem{rahimian2019distributionally}
	Hamed Rahimian and Sanjay Mehrotra.
	\newblock Distributionally robust optimization: A review.
	\newblock {\em arXiv preprint arXiv:1908.05659}, 2019.
	
	\bibitem{gao2017distributional}
	Rui Gao, Xi~Chen, and Anton~J Kleywegt.
	\newblock Distributional robustness and regularization in statistical learning.
	\newblock {\em arXiv preprint arXiv:1712.06050}, 2017.
	
	\bibitem{blanchet2019quantifying}
	Jose Blanchet and Karthyek Murthy.
	\newblock Quantifying distributional model risk via optimal transport.
	\newblock {\em Mathematics of Operations Research}, 44(2):565--600, 2019.
	
	\bibitem{zhao2014data}
	Chaoyue Zhao.
	\newblock {\em Data-driven risk-averse stochastic program and renewable energy
		integration}.
	\newblock University of Florida, 2014.
	
	\bibitem{shafieezadeh2015distributionally}
	Soroosh Shafieezadeh-Abadeh, Peyman~Mohajerin Esfahani, and Daniel Kuhn.
	\newblock Distributionally robust logistic regression.
	\newblock {\em arXiv preprint arXiv:1509.09259}, 2015.
	
	\bibitem{wozabal2014robustifying}
	David Wozabal.
	\newblock Robustifying convex risk measures for linear portfolios: A
	nonparametric approach.
	\newblock {\em Operations Research}, 62(6):1302--1315, 2014.
	
	\bibitem{nguyen2020distributionallyEst}
	Viet~Anh Nguyen, Fan Zhang, Jose Blanchet, Erick Delage, and Yinyu Ye.
	\newblock Distributionally robust local non-parametric conditional estimation.
	\newblock {\em arXiv preprint arXiv:2010.05373}, 2020.
	
	\bibitem{nguyen2020distributionallyMLE}
	Viet~Anh Nguyen, Xuhui Zhang, Jose Blanchet, and Angelos Georghiou.
	\newblock Distributionally robust parametric maximum likelihood estimation.
	\newblock {\em arXiv preprint arXiv:2010.05321}, 2020.
	
	\bibitem{nguyen2021robustifying}
	Viet~Anh Nguyen, Fan Zhang, Jose Blanchet, Erick Delage, and Yinyu Ye.
	\newblock Robustifying conditional portfolio decisions via optimal transport.
	\newblock {\em arXiv preprint arXiv:2103.16451}, 2021.
	
	\bibitem{dou2019distributionally}
	Xialiang Dou and Mihai Anitescu.
	\newblock Distributionally robust optimization with correlated data from vector
	autoregressive processes.
	\newblock {\em Operations Research Letters}, 47(4):294--299, 2019.
	
	\bibitem{qi2021distributionally}
	Meng Qi, Ying Cao, and Zuo-Jun Shen.
	\newblock Distributionally robust conditional quantile prediction with fixed
	design.
	\newblock {\em Management Science}, 2021.
	
	\bibitem{bertsimas2017bootstrap}
	Dimitris Bertsimas and Bart Van~Parys.
	\newblock Bootstrap robust prescriptive analytics.
	\newblock {\em arXiv preprint arXiv:1711.09974}, 2017.
	
	\bibitem{kannan2020residuals}
	Rohit Kannan, G{\"u}zin Bayraksan, and James~R Luedtke.
	\newblock Residuals-based distributionally robust optimization with covariate
	information.
	\newblock {\em arXiv preprint arXiv:2012.01088}, 2020.
	
	\bibitem{gyorfi2006distribution}
	L{\'a}szl{\'o} Gy{\"o}rfi, Michael Kohler, Adam Krzyzak, and Harro Walk.
	\newblock {\em A distribution-free theory of nonparametric regression}.
	\newblock Springer Science \& Business Media, 2006.
	
	\bibitem{nadaraya1964estimating}
	Elizbar~A Nadaraya.
	\newblock On estimating regression.
	\newblock {\em Theory of Probability \& Its Applications}, 9(1):141--142, 1964.
	
	\bibitem{watson1964smooth}
	Geoffrey~S Watson.
	\newblock Smooth regression analysis.
	\newblock {\em Sankhy{\=a}: The Indian Journal of Statistics, Series A}, pages
	359--372, 1964.
	
	\bibitem{tsybakov2008introduction}
	Alexandre~B Tsybakov.
	\newblock {\em Introduction to nonparametric estimation}.
	\newblock Springer Science \& Business Media, 2008.
	
	\bibitem{hannah2010nonparametric}
	Lauren Hannah, Warren Powell, and David Blei.
	\newblock Nonparametric density estimation for stochastic optimization with an
	observable state variable.
	\newblock {\em Advances in Neural Information Processing Systems}, 23:820--828,
	2010.
	
	\bibitem{hanasusanto2013robust}
	Grani~Adiwena Hanasusanto and Daniel Kuhn.
	\newblock Robust data-driven dynamic programming.
	\newblock {\em Advances in Neural Information Processing Systems}, 26:827--835,
	2013.
	
	\bibitem{ban2019big}
	Gah-Yi Ban and Cynthia Rudin.
	\newblock The big data newsvendor: Practical insights from machine learning.
	\newblock {\em Operations Research}, 67(1):90--108, 2019.
	
	\bibitem{ho2019data}
	Chin~Pang Ho and Grani~A Hanasusanto.
	\newblock On data-driven prescriptive analytics with side information: A
	regularized nadaraya-watson approach.
	\newblock {\em URL: http://www. optimization-online. org/DB FILE/2019/01/7043.
		pdf}, 2019.
	
	\bibitem{bolley2007quantitative}
	Fran{\c{c}}ois Bolley, Arnaud Guillin, and C{\'e}dric Villani.
	\newblock Quantitative concentration inequalities for empirical measures on
	non-compact spaces.
	\newblock {\em Probability Theory and Related Fields}, 137(3-4):541--593, 2007.
	
	\bibitem{fournier2015rate}
	Nicolas Fournier and Arnaud Guillin.
	\newblock On the rate of convergence in wasserstein distance of the empirical
	measure.
	\newblock {\em Probability Theory and Related Fields}, 162(3):707--738, 2015.
	
	\bibitem{esteban2020distributionally}
	Adri{\'a}n Esteban-P{\'e}rez and Juan~M Morales.
	\newblock Distributionally robust stochastic programs with side information
	based on trimmings.
	\newblock {\em arXiv preprint arXiv:2009.10592}, 2020.
	
	\bibitem{van2020data}
	Bart~PG Van~Parys, Peyman~Mohajerin Esfahani, and Daniel Kuhn.
	\newblock From data to decisions: Distributionally robust optimization is
	optimal.
	\newblock {\em Management Science}, 2020.
	
	\bibitem{rockafellar2002conditional}
	R~Tyrrell Rockafellar and Stanislav Uryasev.
	\newblock Conditional value-at-risk for general loss distributions.
	\newblock {\em Journal of Banking \& Finance}, 26(7):1443--1471, 2002.
	
	\bibitem{demiguel2007optimal}
	Victor DeMiguel, Lorenzo Garlappi, and Raman Uppal.
	\newblock Optimal versus naive diversification: How inefficient is the 1/n
	portfolio strategy?
	\newblock {\em The Review of Financial Studies}, 22(5):1915--1953, 2007.
	
	\bibitem{fama1992cross}
	Eugene~F Fama and Kenneth~R French.
	\newblock The cross-section of expected stock returns.
	\newblock {\em the Journal of Finance}, 47(2):427--465, 1992.
	
	\bibitem{li2020dimension}
	Wenhao Li, Ningyuan Chen, and L~Jeff Hong.
	\newblock Dimension reduction in contextual online learning via nonparametric
	variable selection.
	\newblock {\em arXiv preprint arXiv:2009.08265}, 2020.
	
	\bibitem{kantorovich1958space}
	Leonid~Vasilevich Kantorovich and SG~Rubinshtein.
	\newblock On a space of totally additive functions.
	\newblock {\em Vestnik of the St. Petersburg University: Mathematics},
	13(7):52--59, 1958.
	
	\bibitem{stein2009real}
	Elias~M Stein and Rami Shakarchi.
	\newblock {\em Real analysis: measure theory, integration, and Hilbert spaces}.
	\newblock Princeton University Press, 2009.
	
	\bibitem{boyd2004convex}
	Stephen Boyd, Stephen~P Boyd, and Lieven Vandenberghe.
	\newblock {\em Convex optimization}.
	\newblock Cambridge University Press, 2004.
	
\end{thebibliography}

\newpage

\newpage
\appendix
\section{Proofs for Section~\ref{sec:theoretical}}
\subsection{Proof of Theorem~\ref{prop:ws-concentration}.}
We first show the result for naive/window kernels.
We use $C_i$, $i=1,2,\dots$ to denote quantities that are independent of $n$, $t$ and $x$. 
Define two auxiliary measures:
    \begin{align*}
        \mu^r_{n,x}&:= \sum_{i=1}^n \frac{ \I{\|x-x_i\|_2\le rh_n}}{\sum_{j=1}^n \I{\|x-x_j\|_2\le rh_n}} \I{y=y_i}\\
        \mu^R_{n,x}&:= \sum_{i=1}^n \frac{ \I{\|x-x_i\|_2\le Rh_n}}{\sum_{j=1}^n \I{\|x-x_j\|_2\le Rh_n}} \I{y=y_i}.
    \end{align*}
    The interpretation of $\mu^r_{n,x}$ and $\mu^R_{n,x}$ are that if we replace $K(\cdot)$ by the naive kernel $\I{\|x\|_2\le r}$ and $\I{\|x\|_2\le R}$, then $\mu_{n,x}$ becomes $\mu^r_{n,x}$ and $\mu^R_{n,x}$. The naive kernels are easier to deal with because they induce empirical distribution over the samples in a ball. By Assumption ~\ref{asp:kernel}, we would try to bound the Wasserstein distance of $\mu_{n,x}$ and $\mu_{Y|x}$ by those of $\mu^r_{n,x}$ ($\mu^{R}_{n,x}$) and $\mu_{Y|x}$.
    
We define another two measures as follows:
    \begin{align*}
        \bar\mu^r_{n,x}&:= \sum_{i=1}^n \frac{ \I{\|x-x_i\|_2\le rh_n}}{\sum_{j=1}^n \I{\|x-x_j\|_2\le rh_n}} \I{Y = y_i-f(x_i)+f(x)}\\
        \bar \mu^R_{n,x}&:=\sum_{i=1}^n \frac{ \I{\|x-x_i\|_2\le Rh_n}}{\sum_{j=1}^n \I{\|x-x_j\|_2\le Rh_n}} \I{Y = y_i-f(x_i)+f(x)}.
    \end{align*}
    The definitions of $\bar\mu^r_{n,x}$ and $\bar\mu^R_{n,x}$ are also intuitive. 
    We are interested in removing the difference between $f(x_i)$ and $f(x)$. Equivalently, the sample points would become $\left\{(x_i, y_i-f(x_i)+f(x))\right\}_{i=1}^n$, and $y_i-f(x_i)+f(x)=\varepsilon_i+f(x)$ are also  i.i.d. samples with the measure $\mu_{Y|x}$. 
    Let $N_n^r$ and $N_n^R$ be the number of data points $x_i$ inside the ball:
    \begin{equation}\label{eq:ball_size}
        N_n^r = \sum_{i=1}^n \I{\|x-x_i\|_2\le rh_n},\quad N_n^R = \sum_{i=1}^n \I{\|x-x_i\|_2\le Rh_n}.
    \end{equation}
    It is obvious that $N_n^r$, $N_n^R$ and $\varepsilon_i$ are independent. Given $N_n^r$ ($N_n^R$), $\bar \mu^r_{n,x}$ ($\bar \mu^R_{n,x}$) is the empirical distribution over $N_n^r$ ($N_n^R$) i.i.d. samples whose covariates are inside the ball. From Theorem 2 in \cite{fournier2015rate}, we have:
    \begin{align}
        \mathbb{P}(\Wscr_p(\mu_{Y|x},\bar \mu^r_{n,x})>t|N_n^r)&\le C_1\exp(-C_2 N_n^r t^2)\label{eq:condition-Nnr}\\
        \mathbb{P}(\Wscr_p(\mu_{Y|x},\bar \mu^R_{n,x})>t|N_n^R)&\le C_1\exp(-C_2 N_n^R t^2)\notag,
    \end{align}
    where $C_1,C_2$ are constants independent of $n$, $x$ and $t$.
    
    Next we would provide a lower probabilistic bound for $N_n^r$ and $N_n^R$. Note that the volume of the intersection of the ball $\|x-x\|_2\le rh_n$ and $[0,1]^{d_x}$ is at least $\min\left\{ 2^{-{d_x}}C_{d_x}\pi r^{d_x} h_n^{d_x}, 2^{-d_x}C_{d_x}\pi\right\}$, where $C_d\pi R^d$ is the volume of a $d$-dimensional ball. Here the worst case is taken when $x$ is in the corner of $[0,1]^{d_x}$. Since the quantity is independent of $t$, $x$ and $n$, we simply use $C_3 h_n^{d_x}$. By Assumption~\ref{asp:distribution}, we have:
    $\mathbb{P}(\|x-x_i\|_2\le rh_n)\ge cC_3 h_n^{d_x}=C_4 h_n^{d_x}$.
    Since $N_n^r$ is the sum of $n$ i.i.d. Bernoulli random variables with success probability $\mathbb{P}(\|x-x_i\|_2\le rh_n)$, we have: $\text{Var}(N_n^r)=\E[N_n^r](1-\E[N_n^r]/n)$.
    From Bernstein's inequality, we would obtain:
    \begin{equation}\label{eq:Nnr}
    \begin{aligned}
        \mathbb{P}(N_n^r\le C_4n h_n^{d_x}/2)&\le \mathbb{P}\left( \frac{1}{n}(N_n^r-\E[N_n^r])\le -\frac{1}{2n}\E[N_n^r]\right)\notag\\
        &\le \exp(- \frac{\E[N_n^r]^2}{4n(2\text{Var}(N_n^r)/n+\E[N_n^r]/3n)})\notag\\
        &\le \exp(- C_5 \E[N_n^r])\le \exp\left(-C_6 nh_n^{d_x}\right)
    \end{aligned}    
    \end{equation}
    Since $N_n^R\ge N_n^r$, we would have the same bound for $N_n^R$. Now combining \eqref{eq:condition-Nnr} and \eqref{eq:Nnr}, we have:
    \begin{align}\label{eq:wm-bar}
        \mathbb{P}(\Wscr_p(\mu_{Y|x},\bar \mu^r_{n,x})>t) &\le \mathbb{P}(\Wscr_p(\mu_{Y|x},\bar \mu^r_{n,x})>t, N_n^r>C_4n h_n^{d_x}/2)+\mathbb{P}(N_n^r\le C_4n h_n^d/2)\notag\\
        &\le C_8\exp(-C_7 nh_n^{d_x} t^2)\notag\\
        \mathbb{P}(\Wscr_p(\mu_{Y|x},\bar \mu^R_{n,x})>t)&\le C_8\exp(-C_7 nh_n^{d_x} t^2)
    \end{align}
    Next we bound $\Wscr_p(\bar \mu^r_{n,x},\mu^r_{n,x})$ and $\Wscr_p(\bar \mu^R_{n,x},\mu^R_{n,x})$. From the definition of $\bar \mu^r_{n,x}$ and $\mu^r_{n,x}$, we would obtain a natural coupling defined as follow:
    \begin{equation*}
        \xi(y_1,y_2) = \sum_{i=1}^n \frac{ \I{\|x-x_i\|_2\le rh_n}}{\sum_{j=1}^n \I{\|x-x_j\|_2\le rh_n}} \I{y_1=y_i-f(x_i)+f(x), y_2=y_i}.
    \end{equation*}
    Therefore,
    \begin{align*}
        \Wscr_p(\bar \mu^r_{n,x},\mu^r_{n,x})&\le \left(\int_{\R \times \R} |y_1-y_2|^p\xi(dy_1,dy_2)\right)^{1/p}\\
                                         &=\left(\sum_{i=1}^n \frac{ \I{\|x-x_i\|_2\le rh_n}}{\sum_{j=1}^n \I{\|x-x_j\|_2\le rh_n}} \|f(x_i)-f(x)\|_2^p\right)^{1/p}\\
                                         &\le rLh_n.
    \end{align*}
    The last inequality follows from Assumption ~\ref{asp:generate}. Similarly, we have:
    \begin{equation*}
        \Wscr_p(\bar \mu^R_{n,x},\mu^R_{n,x})\le RLh_n.
    \end{equation*}
    Combining the above results, we have shown that:
    \begin{align*}
        \mathbb{P}(\Wscr(\mu_{Y|x},\mu^r_{n,x})>t) &\le  \mathbb{P}(\Wscr(\mu_{Y|x},\bar\mu^r_{n,x})>t/2)+ \mathbb{P}(\Wscr(\mu^r_{n,x},\bar\mu_{n,x})>t/2)\\
                                            &\le c_1\exp(-c_2 nh_n^{d_x} t^2)+\I{2rLh_n\le t}.
    \end{align*}
    Therefore we have obtained the result for the naive kernel with radius $r$. Similar arguments hold for the naive kernel with radius $R$. $\hfill\square$

    Based on the two measures $\mu^r_{n,x}$ and $\mu^R_{n,x} (r < R)$, we define the measure as follows:
    \begin{equation*}
    \mu^{R, r}_{n,x}(\lambda) := \sum_{i=1}^n \frac{\lambda \I{\|  x-  x_i\|_2\le rh_n} + (1 - \lambda) \I{\|  x-  x_i\|_2\le Rh_n}}{\sum_{j=1}^n [\lambda \I{\|  x-  x_j\|_2\le rh_n} + (1 - \lambda) \I{\|  x-  x_j\|_2\le Rh_n}]} \I{Y = y_i}.
    \end{equation*}
    which can be interpreted as the linear interpolation of the two naive kernels.
    In order to prove the result for general kernels satisfying Assumption~\ref{asp:kernel}, we first show the following lemma:
\begin{lemma}\label{lemma:ws-concentration-naive-kernel2}
\normalfont
The measure $\mu^{R, r}_{n,x}(\lambda)$ defined above satisfies
    $$\Wscr_p( \mu_{Y | x}, \mu^{R, r}_{n, x}(\lambda)) \leq \max\left\{ \Wscr_p( \mu_{Y | x},\mu^r_{n,x}), \Wscr_p( \mu_{Y | x},\mu^R_{n,x})\right\}, \forall \lambda \in [0, 1].$$
\end{lemma}
\textit{Proof of Lemma \ref{lemma:ws-concentration-naive-kernel2}. } We follow the notation in \eqref{eq:ball_size}. As $r < R$, it is clear that $N_n^r < N_n^R <n$.
We can thus simplify the expression 
\begin{equation*}
    \mu^{R, r}_{n,x}(\lambda) = \sum_{i = 1}^{N_n^r}\frac{1}{N_n^r + (1 -\lambda)(N_n^R - N_n^r)}\I{Y = y_i} + \sum_{i = N_n^r + 1}^{N_n^R}\frac{1 - \lambda}{N_n^r + (1 -\lambda)(N_n^R - N_n^r)}\I{Y = y_i},
\end{equation*}
By definition, we also have $\mu^{R}_{n, x} = \mu^{R, r}_{n,x}(0), \mu^{r}_{n, x} = \mu^{R, r}_{n,x}(1)$. 
We further define:
\begin{equation*}
   \beta_j(\lambda) = \frac{1}{N_n^r + (1 -\lambda)(N_n^R - N_n^r)}\I{j \le N_n^r} + \frac{1 - \lambda}{N_n^r + (1 -\lambda)(N_n^R - N_n^r)}\I{j > N_n^r}, \forall j \in [N_n^R],
\end{equation*}
which satisfies $\sum_{j = 1}^{N_n^R}\beta_j(\lambda) = 1$.

Denote $\gamma \coloneqq \frac{\frac{1}{M_r + (1 - \lambda)(M_R - M_r)} - \frac{1}{M_R}}{\frac{1}{M_r} - \frac{1}{M_R}}$, we then construct a one-to-one mapping from $\lambda \in [0, 1]$ to $\gamma \in [0, 1]$. Then we obtain $\tilde{\beta}_j(\gamma) = (\frac{\gamma}{N_n^r} + \frac{1-\gamma}{N_n^R})\I{j \le N_n^r} + \frac{1-\gamma}{N_n^R}\I{j > N_n^r}, \forall j \in [N_n^R]$, where $\tilde{\beta}_j(\gamma)$ is the weight of $y_i$ for the measure $\mu_{n, x}^{R, r}(\cdot)$.

We first denote the dual representation of Wasserstein distance by Kantorovich's duality in \cite{kantorovich1958space}:
\begin{equation}\label{eq:generaldualwasserstein}
    \Wscr_p^p(\mu, \nu) = \sup_{u, v \in L^1(\Xi)}\{\int_{\Xi}u(\xi)\mu(d\xi) + \int_{\Xi}v(y)\nu(d y): u(\xi) + v(y)\leq \|\xi - y\|_p^p, \forall \xi, y \in \Xi\}.
\end{equation}
For any $\gamma \in [0, 1]$,
\begin{equation*}
    \begin{aligned}
     \Wscr_p^p(\mu_{Y | x}, \mu_{n, x}^{R, r}(\gamma)) =&~\sup_{u, v \in L^1(\Xi)}\sum_{j \in [N_n^R]}\tilde{\beta}_j(\gamma) v(y_i) + \int_{\Xi}u(\xi)\mu(d\xi)\\
   &\text{s.t.}~u(\xi) + v(y_i) \leq \|\xi - y_i\|_p^p, \forall \xi \in \Xi, i \in [n].    
    \end{aligned}
\end{equation*}
We denote the objective functional of the above optimization problem as $g(\gamma; u, v)$, which is linear in $\gamma$, as the feasible region of the function variables $u, v$ is independent of $\gamma$.
Thus, we have the following inequality:
\begin{equation*}
\max_{u, v \in L^1(\Xi)}g(\lambda; {u,v}) \leq \max\{\max_{u, v \in L^1(\Xi)}g(0; {u,v}), \max_{u, v \in L^1(\Xi)}g(1; {u,v})\}.
\end{equation*}
This implies that Lemma~\ref{lemma:ws-concentration-naive-kernel2} holds. $\hfill\square$

We refer to the inequality in Lemma~\ref{lemma:ws-concentration-naive-kernel2} as the \emph{quasi-convex inequality for the naive kernel}. 
As a result, we have:
$$\mathbb{P}(\Wscr_p( \mu_{Y | x},\mu_{n,x}^{R, r}(\lambda))>t)\le \sum_{l \in \{r, R\}}\mathbb{P}(\Wscr_p( \mu_{Y | x},\mu_{n,x}^{l}))>t) \le \tilde{c}_1\exp(-\tilde{c}_2 nt^2), $$
where $\tilde{c}_1, \tilde{c}_2$ are independent with $n, t$ and $  x$.

Next we define the following finite-step kernel as follows:
\begin{equation}
    K(x) = \sum_{i = 1}^{M}\lambda_i\I{\|  x\|_2\le r_i}
\end{equation}
for $r \leq r_1 \leq ...\leq r_M \leq R, \lambda_1 \geq b_r, \sum_{i = 1}^M \lambda_i \leq b_R$, and some $M>0$.
We can prove an analog of Lemma~\ref{lemma:ws-concentration-naive-kernel2} for the finite-step kernel by using the technique in the proof of Lemma~\ref{lemma:ws-concentration-naive-kernel2} recursively.

The finite-step kernel provides an approximation for any kernel between the two naive kernels with width $r$ and $R$. We only show the case of finite-step kernel here.

For any kernel satisfying Assumption~\ref{asp:kernel}, with the approximation theorem for any measure function (i.e. Theorem 4.3) in \cite{stein2009real}, there exists a sequence of step functions converging pointwise almost everywhere. Therefore, for any $K(x)$, there exists a sequence of finite-step kernel functions $\{K_m(x)\}_{m = 1}^{+\infty}$, $K_m(x) \to K(x), \text{a.s.}$. Therefore, for each measure $\mu_{n, x}$ associated with kernel $K(x)$ satisfying Assumption~\ref{asp:kernel}, we can always find a sequence of measure $\{\mu_{n, x}^m\}_{m = 1}^{+\infty}$ associated with finite-step kernel $K_m(x)$, such that $\Wscr_p(\mu_{n, x}, \mu_{n, x}^m) \to 0$ with $m \to \infty$. Combining with the measure concentration result for the finite-step kernel, we would obtain:
\begin{equation}
    \Wscr_p(\mu_{n, x}, \mu_{Y | x}) \leq \Wscr_p(\mu_{n, x}, \mu_{n, x}^m) + \Wscr_p(\mu_{n, x}^m, \mu_{Y | x}).
\end{equation}
Therefore, we provide the measure concentration result for general kernels satisfying Assumption~\ref{asp:kernel}. And we would like to extend the rigorous proof for the general kernel in Assumption ~\ref{asp:kernel} in our future work.

\subsection{Proof of Proposition~\ref{prop:sampleguarantee}.}Under Assumption~\ref{asp:generate}, \ref{asp:distribution} and \ref{asp:kernel}, when $n$ is large enough, if we choose $\varepsilon_n(\alpha):=\frac{\log(\frac{c_1}{\alpha})}{c_2 n}$, we would obtain the following condition that $\mathbb{P}[  \mu_{Y | x} \in B_{\varepsilon}^p(\hat{\mu}_Y)] \geq 1- \alpha$. Therefore, with at least probability $ 1- \alpha$, we would have: $\mathbb{E}_{  \mu_{Y | x}}[c(\hat{z}_n,Y)]\leq \sup_{\mu_{Y} \in B_{\varepsilon}^p(\alpha_n)}\mathbb{E}_{\mu_Y}[c(\hat{z}_n, Y)] = \hat{J}_n$. $\hfill\square$

\subsection{Proof of Proposition~\ref{prop:asymptotic_convergence}.} 
The proof of this proposition is the same as Theorem 3.6 in \cite{esfahani2018data}.

\subsection{Proof of Proposition~\ref{prop:convergence_bound}.} 
We first illustrate a simple form of Wasserstein dual representation when $p = 1$. In \cite{kantorovich1958space}, they show that the type-1 Wasserstein distance between $\mu$ and $\nu$ supported on $\Xi$ admits the following dual representation:
\begin{equation}\label{eq:dualwasserstein}
    \Wscr_1(\mu, \nu) = \sup_{Lip(\phi) \leq 1}\left\{\int_{\Xi} \phi(\xi) \mu(d\xi) - \int_{\Xi} \phi(\zeta) \mu(d\zeta)\right\}.
\end{equation}
When the sample size $n$ is large enough, we obtain the following inequality:
\begin{equation*}
    J^* \leq \E_{Y | x}[c(\hat{z}_n, Y)] \leq \hat{J}_n,a.s.
\end{equation*}
Then we only need to bound the value $\hat{J}_n - J^*$. From the definition of $\Wscr_p(\cdot, \cdot)$, the following inequality holds:
\begin{equation}
    \Wscr_p(\mu_Y, \hat{\mu}_Y) \geq \Wscr_q(\mu_Y, \hat{\mu}_Y), \forall p \geq q, \forall \mu_Y \in \mu(\mathcal{Y}).
\end{equation}
Denote $t:=\|x - y\|$, then Holder inequality implies that $(\E|t|^q)^{\frac{1}{q}} \leq (\E|t|^p)^{\frac{1}{p}}$. Specifically, if we let $p > 1, q = 1$,
we would obtain the following inequality for the cost function:
\begin{equation}
\begin{aligned}
    \hat{J}_n = \inf_{z \in \mathcal{Z}}\sup_{\Wscr_p(\mu_Y, \hat{\mu}_Y) \leq \varepsilon_n(\alpha_n)}\mathbb{E}_{\mu_Y}[c(\hat{z}_n, Y)] &\leq \sup_{\Wscr_p(\mu_Y, \hat{\mu}_Y) \leq \varepsilon_n(\alpha_n)}\mathbb{E}_{\mu_Y}[c(z^*, Y)]\\ &\leq \sup_{\Wscr_1(\mu_Y, \hat{\mu}_Y) \leq \varepsilon_n(\alpha_n)}\mathbb{E}_{\mu_Y}[c(z^*, Y)]\\
    &\leq \sup_{\Wscr_1(\mu_Y, \mu_{Y|x}) \leq 2\varepsilon_n(\alpha_n)}\mathbb{E}_{\mu_Y}[c(z^*, Y)]\\
    &\leq 
    \mathbb{E}_{  \mu_{Y | x}}[c(z^*, Y)] + 2L_z\varepsilon_n(\alpha_n),
\end{aligned}
\end{equation}
where $z^*$ is optimal solution under the true distribution. The first and second inequality holds trivially. The third is given by the fact that the true measure $\mu_{Y |x}$ would stay in the ambiguity set $B_{\varepsilon_n(\alpha_n)}^1(\hat{\mu}_Y)$ when $n$ is large enough. And the fourth inequality follows the dual representation with respect to Wasserstein distance in \eqref{eq:dualwasserstein},i.e.
\begin{equation*}
\begin{aligned}
    |\E_{\mu_Y} [c(z^*, Y)] - \E_{\mu_{Y | x}}[c(z^*, Y)]| &\leq L_z|\E_{\mu_Y} [c(z^*, Y)/L_z] - \E_{\mu_{Y | x}}[c(z^*, Y)/L_z]|\\
    &\leq L_z \Wscr_1(\mu_Y, \mu_{Y | x}) \leq 2L_z\varepsilon_n(\alpha_n)     
\end{aligned}
\end{equation*}
when the measure $\mu_Y$ satisfies $\Wscr_1(\mu_Y, \mu_{Y |x}) \leq 2\varepsilon_n(\alpha_n)$. $\hfill \square$

\section{Proofs for Section~\ref{sec:reformulation}}
\subsection{Proof of Lemma~\ref{lemma:duality}.} 
For general case under this problem structure, the weak duality holds:
\begin{equation*}
    \begin{aligned}
       c_P(z)&= \sup_{\mu_Y \in \mu(\mathcal{Y})}\inf_{\lambda \ge 0}\left\{\E_{\mu_Y}[c(z, y)] + \lambda(\varepsilon^p - \Wscr_p^p(\hat{\mu}_Y, \mu_Y))\right\}\\
       &\le \inf_{\lambda \ge 0}\left\{\lambda \varepsilon^p + \sup_{\mu_Y \in \mu(\mathcal{Y})}\{\E_{\mu_{Y}}[c(z, y)] - \lambda \Wscr_p^p(\hat{\mu}_Y, \mu_Y)\}\right\}.\\
    \end{aligned}
\end{equation*}
The second inequality holds due to the min-max inequality under general cases. For the second item in RHS, based on the duality form of Wasserstein distance representation, we have the following equality:
\begin{equation*}
    \begin{aligned}
       &\sup_{\mu_Y \in \mu(\mathcal{Y})}\left\{\E_{\mu_{Y}}[c(z, y)] - \lambda \Wscr_p^p(\hat{\mu}_Y, \mu_Y)\right\}\\
       =&\sup_{\mu_Y \in \mu(\mathcal{Y})}\left\{\E_{\mu_{Y}}[c(z, y)] - \lambda \sup_{g,  h \in L^1(\Zscr, \Yscr)}\{ \E_{{\mu}_Y}[g(z,y)] + \E_{\xi\sim \hat{\mu}_Y}[h(z,\xi)]\}\right\},\\
    \end{aligned}
\end{equation*}
where $h(z, \xi) \le \inf_{y \in \mathcal{Y}}[\|y- \xi\|_p^p - g(z,y)]$ by \eqref{eq:generaldualwasserstein}. If $\lambda > 0$, we denote  $c_\lambda(\cdot, \cdot): = \frac{c(\cdot, \cdot)}{\lambda} \in L^1(\Zscr, \Yscr)$ and plug $c_\lambda$ into the equation as the inner supremum for $g$ above, we have: 
\begin{equation*}
    \sup_{\mu_Y \in \mu(\mathcal{Y})}\left\{\E_{Y}[c(z, y)] - \lambda \Wscr_p^p(\hat{\mu}_Y, \mu_Y)\right\} \le - \E_{\xi \sim {\mu_Y}, y \sim {\hat{\mu}_Y}}\left\{\inf_{y \in \mathcal{Y}}[\lambda \|\xi - y\|_p^p- c(z, \xi)]\right\}.
\end{equation*}
Therefore the inequality holds for $\lambda > 0$. And for the case of $\lambda = 0$ the inequality trivially holds $c_P(z) \leq c_D(z)$. 

For the strong duality, if we take $\mu_{n, x}(Y)$ to be the nominal distribution $\hat{\mu}_Y$ under Wasserstein ambiguity, then the strong duality holds. The lemma directly refers to Corollary 2 in \cite{gao2016distributionally}, only replace the weight $\frac{1}{n}$ with general weight $\{w_i\}_{i \in [n]} \in \{w_i: w_i \geq 0, \sum_{i \in [n]} w_i = 1\}$. And $\forall y_1, y_2 \in \mathcal{Y}, d^p(y_1, y_2) = \|y_1 - y_2\|_p^p$ in their corollary. Therefore our lemma holds.$\hfill \square$

\subsection{Proof of Proposition~\ref{prop:regularization}.}This proposition holds generally for $d_y \in \Zscr^+$. We make the following regularity assumptions for $c(z, y)$ similar to \cite{gao2017distributional}:
For completeness, we list the assumptions in \cite{gao2017distributional} in our notation.
\begin{assumption}\label{asp:regular_upper}
\normalfont
     $c(z, \cdot)$ is Lipschitz continuous in $z$.
\end{assumption}
\begin{assumption}\label{asp:regular_lower}
    $c(z, \cdot)$ is differentiable and there exists constant $\kappa \in [0, 1]$ and $h_z(y): \mathcal{Y} \to \mathbb{R}$ is uniformly bounded in $\mathcal{Y}$ such that:
    \begin{equation}
        |\nabla_y c(z, y_1) - \nabla_y c(z, y_2)| \leq h_z(y_2) | y_1 - y_2 |^{\kappa}, \forall y_1, y_2 \in \mathcal{Y}.  
    \end{equation}    
\end{assumption}

From Assumption \ref{asp:regular_upper} and proofs in Proposition~\ref{prop:convergence_bound}, for any given $z$, the following conclusion holds:
\begin{equation}\label{eq:upper}
    \sup_{\Wscr_p(\mu_Y, \hat{\mu}_Y) < \varepsilon_n(\alpha_n)}\mathbb{E}_{\mu_Y}[c(z, Y)] \leq \mathbb{E}_{\hat{\mu}_Y}[c(z, Y)] + 2L_z\varepsilon_n(\alpha_n).
\end{equation}
This presents a upper bound for the worst-case cost.

On the other hand, by Assumption \ref{asp:regular_lower}, we claim the following two lower bounds for the worst-case cost classified by Wasserstein order $p$:

If $p \in (\kappa + 1, \infty]$, we have:
\begin{equation}\label{eq:lower1}
    \sup_{\Wscr_p(\mu_Y, \hat{\mu}_Y) < \varepsilon_n(\alpha_n)}\mathbb{E}_{\mu_Y}[c(z, Y)] \geq 
    \mathbb{E}_{\hat{\mu}_Y}[c(z, Y)] + \alpha_n \| \nabla_y c(z, y)\|_{\hat{\mu}_Y, p_*} - \alpha_n^{\kappa + 1}\|h_z(y)\|_{\hat{\mu}_Y, \frac{p}{p - \kappa - 1}}.
\end{equation}
If $p \in [1, \kappa + 1]$, we have:
\begin{equation}\label{eq:lower2}
    \sup_{\Wscr_p(\mu_Y, \hat{\mu}_Y) < \varepsilon_n(\alpha_n)}\mathbb{E}_{\mu_Y}[c(z, Y)] \geq 
    \mathbb{E}_{\hat{\mu}_Y}[c(z, Y)] + \alpha_n \| \nabla_y c(z, y)\|_{\hat{\mu}_Y, p_*} - \alpha_n^{\kappa + 1}\|h_z(y)\|_{\hat{\mu}_Y, \infty}. 
\end{equation}
where the empirical norm of a function in \eqref{eq:lower1} and \eqref{eq:lower2} is defined below:
\begin{equation}
    \| h_z(y) \|_{\hat{\mu}_Y, p}:= 
    \begin{cases}
    (\sum_{i \in [n]}w_i\|h_z(y_i)\|_2^{p})^{\frac{1}{p}},&1\leq p < +\infty, \\
    \max_{i \in [n]}\|h_z(y_i)\|_2,&p = \infty.
    \end{cases}
\end{equation}
The proof of \eqref{eq:lower1} and \eqref{eq:lower2} is similar to Proposition 2 in \cite{gao2017distributional}. 
To be more specific, considering any distribution supported on $n$ discrete points $\{\xi_i\}_{i \in [n]}$, we state that $\sup_{\Wscr_p(\mu_Y, \hat{\mu}_Y) < \varepsilon_n(\alpha_n)}\mathbb{E}_{\mu_Y}[c(z, Y)]$ is lower bounded by the optimal value of the following objective function following \eqref{eq:duality2} in Lemma \ref{lemma:duality}:
\begin{equation}\label{eq:regular_mid_target}
\sup_{\xi_i \in \mathcal{Y}}\left\{\sum_{i = 1}^{n}w_i [c(z, \xi_i) - c(z, y_i)]: (\sum_{i = 1}^{n}w_i\|\xi_i- y_i\|_p^p)^{1/p} \leq \varepsilon\right\},    
\end{equation}
If the cost function $c(z,y)$ satisfies Assumption \ref{asp:cost_function}, then the lower bound above can be replaced with the equality. Furthermore, from $\nabla_y c(z, y)$ in Assumption \eqref{asp:regular_lower} and mean-value theorem, the value in \eqref{eq:regular_mid_target} can be lower bounded as follows:
\begin{equation}\label{eq:regular_mid_step}
\begin{aligned}
&\sup_{\xi_i \in \mathcal{Y}}\left\{\sum_{i = 1}^{n}w_i [ \|\nabla_y c(z, y_i)\|_p \|\xi_i - y_i\|_p - h_z(y_i)\|\xi_i - y_i\|_p^{\kappa + 1}]: (\sum_{i = 1}^{n}w_i\|\xi_i- y_i\|_p^p)^{1/p} \leq \varepsilon\right\}\\
=&\sup_{\xi_i \in \mathcal{Y}}\left\{\sum_{i = 1}^{n}w_i \|\nabla_y c(z, y_i)\|_p \|\xi_i - y_i\|_p: (\sum_{i = 1}^{n}w_i\|\xi_i- y_i\|_p^p)^{1/p} \leq \varepsilon\right\}\\
-&\sup_{\xi_i \in \mathcal{Y}}\left\{\sum_{i = 1}^{n}h_z(y_i)\|\xi_i - y_i\|_p^{\kappa + 1}: (\sum_{i = 1}^{n}w_i\|\xi_i- y_i\|_p^p)^{1/p} \leq \varepsilon \right\}\\
\end{aligned}   
\end{equation}
For the first item in RHS of \eqref{eq:regular_mid_step},  Holder inequality shows that:
\begin{equation}
\begin{aligned}
\sum_{i = 1}^{n}w_i \|\nabla_y c(z, y_i)\|_p \|\xi_i - y_i\|_p &\leq [\sum_{i = 1}^{n} (w_i^{1/{p^*}} \|\nabla_y c(z, y_i)\|_p)^{p^*}]^{1/{p^*}}[\sum_{i = 1}^{n} (w_i^{1/p} \|\xi_i - y_i\|_p)^p]^{1/p} \\
&\leq \|\nabla_y c(z, y)\|_{\hat{\mu}_Y, p^*}\alpha_n.
\end{aligned}
\end{equation}
And the second item of RHS above can be bounded from the definition of the empirical norm:
\begin{equation}
    RHS_2\leq 
    \begin{cases}
    \alpha_n^{\kappa + 1}\|h_z(y)\|_{\hat{\mu}_Y, \frac{p}{p - \kappa - 1}},&p > \kappa + 1 \\
    \alpha_n^{\kappa + 1}\|h_z(y)\|_{\hat{\mu}_Y, \infty},&p \leq \kappa + 1.
    \end{cases}
\end{equation}
Therefore, we provide two cases for the worst-case cost with respect to $p$ in \eqref{eq:lower1} and \eqref{eq:lower2}. Then with the same method applied to Theorem 2 in \cite{gao2017distributional}, Borel-Cantelli lemma can show that $\lim_{n \to \infty}\|\nabla_y c(z, y)\|_{\hat{\mu}_Y, p^*} = L_z$, a.s. The convergence actually holds for all $p$. Therefore combining \eqref{eq:upper}, \eqref{eq:lower1} and \eqref{eq:lower2}, we then prove the proposition here.$\hfill\square$
\subsection{Proof of Proposition~\ref{prop:reformulation}.} Convexity holds under Assumption \ref{asp:cost_function_piecewise}. Based on Lemma \ref{lemma:duality} and Wasserstein order $p = 1$, $c_p(z) = \min_{\lambda \geq 0}\{\lambda \varepsilon + \sum_{i \in [n]} w_i\sup_{y \in \mathcal{Y}}[c(z, y) - \lambda |y - y_i|]\}$

Introducing the auxiliary variable $s_i \geq \sup_{y \in \mathcal{Y}}[c(z, y) - \lambda |y - y_i|]$, then the objective function is linear in $\lambda, \{s_i\}_{i \in [n]}$ but has infinite constraints. 
We deal with $s_i$ as follows:
\begin{equation}
\begin{aligned}
s_i &\geq \sup_{y \in \mathcal{Y}}[c(z, y) - \lambda |y - y_i|\\
&=\sup_{y \in \mathcal{Y}}[c_k(z, y) - \max_{|p_{ik}| \leq \lambda }p_{ik}(y - y_i)], \forall i \in [n], k \in [K].\\
&=\min_{|p_{ik}| \leq \lambda }\max_{y \in \mathcal{Y}}[c_k(z, y) - p_{ik}(y - y_i)], \forall i \in [n], k \in [K],\\
\end{aligned}
\end{equation}
where the second equality holds due to the Cauchy inequality and the fact that $c(z, Y) = \max_{k \in [K]} c_k(z, Y)$. The third equality holds because of Sion's minimiax theorem for convex functions.

From the definition of convex conjugate, it is clear that:
\begin{equation*}
\begin{aligned}
    \sup_{y \in \mathcal{Y}}[c_k(z, y) - p_{ik}(y - y_i)] &= [-c_{z, k} + \I{\cdot \in \mathcal{Y}}]^*(p_{ik})\\
    &=\inf_{p_{ik}}{[-c_{z, k}]^*(p_{ik} - v_{ik}) + \sigma_{\mathcal{Y}}^*(v_{ik})}\\
\end{aligned}
\end{equation*}
The first equality holds because of conjugate of involution property in Chapter 3 in \cite{boyd2004convex}. The second equality holds because $c_k(z, Y)$ is upper semi-continuous. Therefore, we reach the conclusion in Proposition \ref{prop:reformulation}. $\hfill\square$

For our numerical studies, we represent this piece-wise linear model more specifically. The when uncertainty parameter $y \in \mathcal{Y} (\subset \mathbb{R}^{d_y}): = \{A y \leq b\}$, where $A \in \mathbb{R}^{d_1\times d_y}, b \in \mathbb{R}^{d_1}$. We also assume $c_{z, k}(y) = f_{z, k} y + g_{z, k}$, where $f_{z, k}, g_{z, k} \in \R$ under the cost function $c(z, Y) = \max_{k \in [K]}c_k(z, Y)$. Then $c_P(z)$ is equivalent to:
\begin{equation}\label{eq:linear_reformulation}
\begin{aligned}
\inf_{\lambda, s_i, \gamma_{ik}}~&\lambda \varepsilon + \sum_{i = 1}^{n}w_is_i\\
\text{s.t.}~&g_{z, k} + f_{z, k}y + \gamma_{ik}^{\top}(b - Ay) \le s_i, \forall i \in [n], \forall k \in [K]\\
~& \| A^{\top}\gamma_{ik} - f_{z, k}\|_{\infty} \le \lambda, \forall i \in [n], k \in [K]\\
~&\gamma_{ik} \geq 0, \forall i \in [n], k \in [K].
\end{aligned}
\end{equation}

\eqref{eq:linear_reformulation} is the direct corollary \ref{prop:reformulation}. For general $d_y \geq 1$, $p_{ik}$ in \eqref{eq:reformulation} need be revised to $\|p_{ik}\|_{\infty} \leq \lambda, \forall i \in [n], k \in [K]$. For the items in the first constraint, strong duality implies that:
\begin{equation*}
\begin{aligned}
    \sigma_{\mathcal{Y}}(v_{ik}) &= \{\sup_{y} v_{ik}^{\top} y,~\text{s.t.}A y \leq b\}\\
    &=\{\inf_{\gamma_{ik} \geq 0} \gamma_{ik}^{\top}b, ~\text{s.t.} A^{\top} \gamma_{ik} = v_{ik}\}.
\end{aligned}
\end{equation*}
And the definition of conjugate function demonstrates that:
\begin{equation*}
[-c_{z, k}]^*(p_{ik} - v_{ik}) = [-f_{z.k}]^*(p_{ik} - v_{ik}) = 
\begin{cases}
    -g_{z, k},& p_{ik} - v_{ik} = -f_{ik}, \\
    \infty,& p_{ik} - v_{ik} \not= -f_{ik}
\end{cases}
\end{equation*}
Then $p_{ik} = A^{\top}\gamma_{ik}$ in . we plug those above into \eqref{eq:reformulation} to obtain \eqref{eq:linear_reformulation}.

We demonstrate that the two problem case in Section~\ref{sec:numerical} can be reformulated into the linear programming problem.

\textit{Newsvendor. }
Here $c(z, y) = \max \{b\cdot y - bz, -h\cdot y + hz\}$, and $A = [-1;1]^{\top}, b = (-L, U)^{\top}$ in \eqref{eq:linear_reformulation}, where the feasible region of demand $Y$ is $\mathcal{Y} = \{Y \in \mathbb{R}: L \leq Y \leq U\}.$ Therefore, the model can be reconstructed below with history samples $\{y_i\}_{i \in [n]}$ as well as kernel weights $\{w_i\}_{i \in [n]}$.
\begin{equation*}
\begin{aligned}
    \min_{z, \lambda, s_i, \gamma_{ik}}&~\lambda \varepsilon + \sum_{i = 1}^{n}w_i s_i\\
    \text{s.t.}&~ - bz + (y_i - L) \gamma_{i1} + (U - y_i)\gamma_{i2} \leq s_i - by_i, \forall i \in [n],\\
    &~hz + (y_i - L) \gamma_{i3} + (U - y_i)\gamma_{i4} \leq s_i + hy_i, \forall i \in [n],\\
    &~ -\lambda + b \leq -\gamma_{i1} + \gamma_{i2} \leq \lambda + b, \forall i \in [n],\\
    &~ -\lambda - h \leq -\gamma_{i3} + \gamma_{i4} \leq \lambda - h, \forall i \in [n],\\
    &~\gamma_{ik} \geq 0, \forall i \in [n], k \in \{1,2,3,4\}, z \in \mathcal{Z}, \lambda \geq 0.
\end{aligned}
\end{equation*}

\textit{Portfolio Allocation. }We note that the decision variable is $z^{\prime}:=(z, v)$ by the min-max rule with the cost function below:
\begin{equation*}
\begin{aligned}
\inf_z \sup_{\mu_y \in B_{\varepsilon}^p(\hat{\mu}_Y)}c(z, y) &= \min_z \sup_{\mu_y \in B_{\varepsilon}^p(\hat{\mu}_Y)} [\text{CVaR}_{\eta}(y^{\top}z) - \gamma y^{\top} z]\\
& = \inf_{z, v}\{v - \gamma y^{\top}z + \frac{1}{\eta}\sup_{\mu_y \in B_{\varepsilon}^p(\hat{\mu}_Y)}\mathbb{E}_{y}[\max\{-y^{\top} z - v, 0\}]\}\\
& = \min_{z, v} \max_{\mu_y \in B_{\varepsilon}^p(\hat{\mu}_Y)}\mathbb{E}_y[\max\{-(\gamma + \frac{1}{\eta} y^{\top}z) + (1 - \frac{1}{\eta}) v, -\gamma y^{\top} z + v\}]\\
\end{aligned}
\end{equation*}
Here we do not impose linear constraints for the distribution of assets' returns. Then $A = 0 \in \mathbb{R}^{d_1 \times d_y}$, $b = 0$, it follows that the problem is reformulated below:
\begin{equation}
\begin{aligned}
    \min_{  z \in \Zscr, v, \eta, s_i}&~\lambda \varepsilon+ \sum_{i = 1}^{n}w_i s_i\\
    \text{s.t.}&~ - (\frac{1}{\eta} + \gamma)  y_i^{\prime}  z + v(1-\frac{1}{\eta}) \le s_i, \forall i \in [n],\\
    &~ -\gamma   y_i^{\prime}   z + v \le s_i, \forall i \in [n],\\
    &~ -\lambda \le (\frac{1}{\eta} + \gamma)z_i\le \lambda, \forall i \in [n],\\
    &~ -\lambda \le \gamma z_i \le \lambda, \forall i \in [n].
\end{aligned}
\end{equation}

\end{document}